\documentclass[]{interact}

\usepackage{epstopdf}
\usepackage[caption=false]{subfig}

\usepackage[numbers,sort&compress]{natbib}
\bibpunct[, ]{[}{]}{,}{n}{,}{,}
\makeatletter
\def\NAT@def@citea{\def\@citea{\NAT@separator}}
\makeatother
\usepackage{anyfontsize}
\theoremstyle{plain}
\newtheorem{theorem}{Theorem}[section]
\newtheorem{lemma}[theorem]{Lemma}
\newtheorem{corollary}[theorem]{Corollary}

\theoremstyle{definition}

\newtheorem{example}[theorem]{Example}

\usepackage{xcolor}

\theoremstyle{remark}

\begin{document}
	

	\title{Spectral Properties of Off-Diagonal Block Linear Relations via Moore–Penrose Inverses in Hilbert Spaces}
  \author{\name {Arup Majumdar \thanks{Arup Majumdar (corresponding author). Email address: arupmajumdar93@gmail.com}} \affil{Institute of Mathematics, Czech Academy of Sciences,
			\v{Z}itn\'a 25, Prague, Czech Republic.}}
  
	\maketitle

	\begin{abstract}
In this paper, we characterize the essential spectra and the resolvent set of the off-diagonal block linear relation
\[
\begin{bmatrix}
0 & A \\
B & 0
\end{bmatrix}
\]
in terms of the essential spectra and resolvent sets of the products $AB$ and $BA$. Our approach establishes precise spectral relationships that connect the structural properties of the block linear relation with those of the associated compositions.

Furthermore, we investigate the Moore--Penrose inverses of closed linear relations in Hilbert spaces and employ these results to extend the spectral analysis to the off-diagonal block linear relation
\[
\mathcal{A} =
\begin{bmatrix}
0 & T^{\dagger} \\
T & 0
\end{bmatrix},
\]
where $T$ is a closed, continuous linear relation with closed range from a Hilbert space $H$ to a Hilbert space $K$, and $T^{\dagger}$ denotes its Moore--Penrose inverse.
\end{abstract}

      \begin{amscode}47A06; 47A10, 47A53.\end{amscode}
    \section{Introduction}
    The Moore-Penrose inverse is a key concept in linear algebra and functional analysis, providing a generalized inverse for matrices and linear operators that may not have a traditional inverse. Named after E. H. Moore and Roger Penrose, who independently introduced the concept in 1920 and 1955, respectively, it plays a crucial role in understanding operator behaviour. Additionally, Yu. Ya. Tseng extends the idea of Moore-Penrose (generalized) inverses to densely defined linear operators on Hilbert spaces in his works \cite{MR0029479, MR0031192, MR0031191}. The study of Moore-Penrose inverses of closed operators is explored in depth in \cite{MR0136996, MR0473881}. In the paper \cite{majumdar2024hyers}, the authors introduce the Moore-Penrose inverses of closable operators with decomposable domains while establishing Lemma 2.10 \cite{majumdar2024hyers}. Moreover, the paper \cite{MR0396607, MR007} provides several characterizations of Moore-Penrose inverses for closable operators in Hilbert spaces.\\
   Recent developments also include the Moore-Penrose inverse of linear relations in Hilbert spaces. In \cite{MR3231505}, the author introduces the Moore-Penrose inverse for closed relations and approximates it by a sequence of bounded, finite-rank operators. A critical point here is that the range of a closed operator must be dense in the Hilbert space to derive properties of its Moore-Penrose inverse. However, there are closed operators whose domains are not dense in the space, which raises the need for characterizations of Moore-Penrose inverses in the context of such operators. This is especially relevant because the adjoint of a linear relation can exist without requiring a dense domain.\\
   The Moore-Penrose inverse has been widely studied due to its applications in areas such as multiple regression, linear systems theory, singular value decomposition, and digital image restoration \cite{greville1958pseudoinverse, chountasis2009applications}. This paper delves into the exploration of several interesting properties of Moore-Penrose inverses of closed linear relations in Hilbert spaces. Section 2 is dedicated to the basic definitions and notations related to linear relations. In Section 3, we discuss the intriguing results that characterize the Moore-Penrose inverses of closed linear relations. In the final section, we characterize the essential spectra and the resolvent set of the off-diagonal block linear relation in terms of the essential spectra and resolvent sets of the compositions of its two off-diagonal entries. In particular, we apply these results to a special class of off-diagonal block linear relations in which one off-diagonal entry is a linear relation and the other is its Moore–Penrose inverse.

\section{Preliminaries}
Throughout the paper, the symbols $H, K, H_{i}, K_{i} ~(i = 1, 2)$ represent real or complex Hilbert spaces. A linear relation $T$ from $H$ into $K$ is a linear subspace of the Cartesian product in $H \times K$, and the collection of all linear relations from $H$ into $K$ is denoted by $LR(H, K)$. We call $T$ a closed linear relation from $H$ into $K$ if it is a closed subspace of $H \times K$ and the set of all closed linear relations from $H$ into $K$ is denoted by $CR(H, K)$. The following notations of domain, range, kernel and multi-valued part of a linear relation $T$ from $H$ into $K$ will be used respectively in the paper:\\
 $$D(T) = \{h \in H: \{h, k\} \in T\},  ~ R(T) = \{k \in K: \{h, k\} \in T\}$$
 $$ N(T)= \{h\in H:\{h, 0\} \in T\}, ~ M(T) = \{k\in K: \{0, k\} \in T \}.$$
 It is obvious that $N(T)$ and $M(T)$ both are closed subspaces in $H$ and $K$, respectively, whenever $T$ is a closed linear relation from $H$ into $K$. In general, the inverse of an operator is always a linear relation. The inverse of a linear relation $T$ from $H$ into $K$ is defined as $T^{-1}= \{\{k, h\} \in K \times H: \{h, k\} \in T \}$. Thus, it is immediate that $D(T^{-1}) = R(T)$ and $N(T^{-1}) = M(T)$. We define $Tx = \{y \in K: \{x, y\} \in T\}, \text{ where } T \in LR(H, K).$ Consequently, $T\vert_{W}$ denotes the restriction of $T \in LR(H, K)$ with domain $D(T) \cap W$, where $W$ is a subset of $H$ (in other words, $T\vert_{W}$ is equal to $T$ in domain $D(T) \cap W$). The adjoint of a linear relation $T$ from $H$ into $K$ is the closed linear relation $T^{*}$ from $K$ into $H$ defined by:
 $$T^{*} = \{\{k^{'}, h^{'}\} \in K \times H : \langle h^{'}, h \rangle = \langle k^{'}, k \rangle, \text{ for all } \{h, k\} \in T\}.$$
 Observe that $(T^{-1})^{*} = (T^{*})^{-1}$, so that $(D(T))^{\perp} = M(T^{*})$ and $N(T^{*}) = (R(T))^{\perp}$. A linear relation $T$ in $H$ is said to be symmetric if $T \subset T^{*}$. Again, a linear relation $T$ in $H$ is non-negative if $\langle k, h \rangle \geq 0$, for all $\{h, k\} \in T$. Moreover, a linear relation $T$ in $H$ is said to be self-adjoint when $T = T^{*}$. If $S$ and $T$ both are linear relations then their product $TS$ defined by:
 $$TS = \{\{x,y\}: \{x,z\} \in S \text{ and } \{z,y\} \in T, \text{ for some } z \}.$$
 The sum of two linear relations $T$ and $S$ from $H$ into $K$ is $$T + S = \{\{x, y +z\}: \{x, y\} \in T \text{ and } \{x, z\} \in S \}.$$ Whereas, the Minkowski sum is denoted by $$T \widehat{+} S =\{\{x+v, y+w\}:\{x, y\} \in T,\{v, w\} \in S\}.$$ Consider $T \in CR(H, H)$, a point $\lambda \in \mathbb{C}$ is said to belong to the resolvent set $\rho(T)$ of $T$ if $(T - \lambda)^{-1}$ is a bounded operator in the domain $H$ and the spectrum $\sigma(T)$ of $T$ is the complement of $\rho(T)$ in $\mathbb{C}.$\\
 Here, $Q_{T}$ denotes the natural quotient map from $K$ into $K / \overline{M(T)}$, where $T$ is a linear relation from $H$ into $K$. It is easy to show that $Q_{T}T$ is a linear operator from $H$ into $K / \overline{M(T)}$. We call the linear relation $T$ from $H$ into $K$ a continuous linear relation if $Q_{T}T$ is a bounded operator, and the set of all continuous linear relations from $H$ into $K$ is denoted by $BR(H, K).$ If $D(T) = H$ and $T \in BR(H, K),$ then we shall say that $T$ is bounded, and the set of all bounded linear relations from $H$ into $K$ is denoted by $BCR(H, K)$. We denote $\mathcal{KR}(H,K)$ the class of compact linear relations from \(H\) to \(K\).
A linear relation \(T \in {LR}(H,K)\) is called compact if
$\overline{Q_T(T(B_H))}$ is compact in \(K\), where \(B_H\) denotes the unit ball of \(H\).\\
A closed linear relation \(T\) from \(H\) into \(K\) is said to be an
upper semi-Fredholm relation, abbreviated by \(\Phi_+\), if \(T\) has closed range and
\( \dim(N(T)) = \alpha(T) < \infty\). In this case, we write \(T \in \Phi_+(H,K)\). The relation \(T\) is called a lower semi-Fredholm relation, abbreviated by \(\Phi_-\), if \(R(T)\) is a closed subspace of \(K\) with finite codimension. We then write \(T \in \Phi_-(H,K)\).
The relation \(T\) is said to be a semi-Fredholm relation (respectively, a Fredholm relation)
if
\[
T \in \Phi_+(H,K) \cup \Phi_-(H,K)
\quad
(\text{respectively, }
T \in \Phi_+(H,K) \cap \Phi_-(H,K)).
\]\\
 The regular part of a closed linear relation $T$ from $H$ into $K$ is $P_{\overline{D(T^{*})}}T$, denoted by $T_{op}$ which is an operator with $T_{op} \subset T$, where $P_{\overline{D(T^{*})}}$ is the orthogonal projection in $K$ onto $\overline{D(T^{*})} = (M(T))^{\perp}$. It can be shown that $T = T_{op} \widehat{+} ~(\{0\} \times M(T))$, when $T$ is a closed linear relation from $H$ into $K$ \cite{MR3971207}. When $T$ is a closed relation from $H$ into $K$, then $T_{op}$ is also a closed operator \cite{MR3971207}. In \cite{MR3971207}, the Moore-Penrose inverse of the linear relation $T$ from $H$ into $K$ is defined by the linear operator $(T^{-1})_{op} = P_{(N(T))^{\perp}}T^{-1}$, where the domain of $(T^{-1})_{op}$ is $R(T)$. However, throughout, we will follow the general definition of the Moore-Penrose inverse of the closed linear relation $T \in CR(H, K)$, which is defined by $P_{(N(T))^{\perp}} T^{-1}P_{R(T)}$ in the domain $R(T) \oplus (R(T))^{\perp}$, is denoted by $T^{\dagger}$ \cite{MR3231505}. In general, $P_{W}$ is the orthogonal projection on the closed subspace $W$ of the specified Hilbert space and $G(P_{W})$ is the graph of the orthogonal projection $P_{W}$. However, we will consider the notation $P_{W}$ as the graph of the orthogonal projection from the specified Hilbert space onto its closed subspace $W$. It is noted that $R(T)$ is closed if and only if $R(T^{*})$ is closed if and only if $R(TT^{*})$ is closed if and only if $R(T^{*}T)$ is closed, when $T \in CR(H, K)$ \cite{MR3079830}. Furthermore, $R(T) = R(TT^{*})$ and $R(T^{*}T) = R(T^{*})$, when $T \in CR(H, K)$ is a  closed range linear relation.

 \section{Characterization of the Moore-Penrose inverses of closed linear relations}
 From now on, $T$ is a closed linear relation from Hilbert space $H$ into $K$. So, the Moore-Penrose inverse of $T$ is denoted by $T^{\dagger} =  P_{(N(T))^{\perp}} T^{-1}P_{R(T)}$ in the domain $R(T) \oplus (R(T))^{\perp}$. The most interesting feature of $T^{\dagger}$ is the element $\{x, y\}$ of $T^{\dagger}$ can be written as $\{x, x^{'}\} \in P_{R(T)}$, $\{y, x^{'}\} \in T$ and $\{y, y\} \in P_{(N(T))^{\perp}}$, where $x = x^{'} + x^{''}$, for some $x^{'} \in R(T)$, $x^{''} \in (R(T))^{\perp}$.
 \begin{theorem}\label{thm 3.1}
 Let $T \in CR(H, K)$ be a closed range linear relation. Then the Moore–Penrose inverse of $T$ satisfies $$(T^{\dagger})^{*} = (T^{*})^{\dagger}.$$
 \end{theorem}
 \begin{proof}
 Since $R(T)$ is closed, it follows that $R(T^{*})$ is closed. Consequently, the Moore–Penrose inverse of $T^{*}$ exists and given by $$(T^{*})^{\dagger} = P_{R(T)}(T^{*})^{-1} P_{R(T^{*})}.$$ By Proposition 1.3.9(ii) \cite{MR3971207}, we compute
 \begin{align*}
 (T^{\dagger})^{*} &= (P_{(N(T))^{\perp}} T^{-1} P_{R(T)})^{*}\\
 &= (T^{-1}P_{R(T)})^{*}P_{(N(T))^{\perp}}~ ( \text{$P_{(N(T))^{\perp}}$ is a bounded operator and $T^{-1}P_{R(T)}$ is closed})\\
 &= (T^{-1}P_{R(T)})^{*}P_{R(T^{*})}.
 \end{align*}
 It is immediate from Proposition 1.3.9 \cite{MR3971207} that $$(T^{-1}P_{R(T)})^{*} \supset P_{R(T)} (T^{-1})^{*}.$$ To establish the reverse inclusion, let $\{h,k\} \in (T^{-1}P_{R(T)})^{*}.$ Then, $$\langle h, x \rangle = \langle k, y\rangle, ~\text{for all}~ \{y, x\} \in T^{-1}P_{R(T)}.$$ For every $y \in R(T)^{\perp}$, we have $\{y, 0\} \in T^{-1}P_{R(T)}$ and $\langle k, y \rangle = 0,$ which implies that $k \in R(T).$ Similarly, for every $x \in N(T)$, we get $\{0,x\} \in T^{-1}P_{R(T)}$ yields $\langle h, x \rangle = 0,$ thus $h \in N(T)^{\perp} = R(T^{*}).$ Now, if $\{u, v\} \in T^{-1}$, then $u \in R(T)$ hence $\{u, v\} \in T^{-1}P_{R(T)}.$ This shows that $$T^{-1} \subset T^{-1}P_{R(T)}$$ and therefore, $$(T^{-1}P_{R(T)})^{*} \subset (T^{-1})^{*}.$$ Consequently $\{h, k\} \in (T^{-1})^{*}$, $\{k, k\} \in P_{R(T)}$ which implies $\{h, k\} \in P_{R(T)}(T^{-1})^{*}$. So, we get the equality $(T^{-1}P_{R(T)})^{*} = P_{R(T)}(T^{-1})^{*}.$ Finally, combining the above identities, we obtain
  \begin{align*}
  (T^{\dagger})^{*} = P_{R(T)}(T^{*})^{-1}P_{R(T^{*})} = (T^{*})^{\dagger}.
  \end{align*}
 \end{proof}
 Theorem \ref{thm 3.2} has been proven in Proposition 3.3 \cite{MR3231505}. However, a different approach is mentioned to prove Theorem \ref{thm 3.2}.
 \begin{theorem}\label{thm 3.2}
 Let $T \in CR(H, K)$. Then the Moore–Penrose inverse $T^{\dagger}$ is a closed operator. Moreover, if $R(T)$ is closed, then $T^{\dagger}$ is bounded.
 \end{theorem}
 \begin{proof}
  Let $\{x, y\} \in \overline{T^{\dagger}}.$ Then there exists a sequence $\{\{x_{n}, y_{n}\}\}$ in $T^{\dagger}$ such that $$\{x_{n}, y_{n}\} \to \{x, y\} ~\text{as}~ n \to \infty,$$ that is, $$x_{n} \to x~ \text{and}~ y_{n} \to y ~ \text{as}~ n \to \infty.$$ For all $n \in \mathbb{N},$ decompose  $$x = x^{'} + x^{''}, ~~~~ x_{n} = x_{n}^{'} + x_{n}^{''},$$ where $$x^{'} \in \overline{R(T)},~ x_{n}^{'} \in R(T), ~\text{and}~ x^{''}, x_{n}^{''} \in (R(T))^{\perp}.$$ Then $x_{n}^{'} \to x^{'}$, as $n \to \infty$, Consequently, $$\{x_{n}, x_{n}^{'}\} \in P_{R(T)}, ~~~ \{x_{n}^{'}, y_{n}\} \in T^{-1}~~~ \{y_{n}, y_{n}\} \in P_{(N(T))^{\perp}}.$$ Since, $T^{-1}$ is closed, passing to the limit yields $$\{x, x^{'} \} \in P_{R(T)},~~ \{x^{'}, y\} \in T^{-1},~~ \{y, y\} \in P_{(N(T))^{\perp}}$$ which implies $\{x, y\} \in T^{\dagger},$ showing that $T^{\dagger}$ is closed on its domain $R(T) \oplus (R(T))^{\perp}$. Next, suppose that $\{0,h\} \in T^{\dagger}.$ Then $$\{0, 0\} \in P_{R(T)},~~ \{0, h\} \in T^{-1},~~ \{h, h\} \in P_{(N(T))^{\perp}}.$$ This implies that $h \in N(T) \cap (N(T))^{\perp} = \{0\}.$ Therefore, $T^{\dagger}$ is single valued and hence a closed operator.\\ Finally, when $R(T)$ is closed, then by using the closed graph Theorem, we obtain that $T^{\dagger}$ is bounded.
 \end{proof}
 \begin{example}
 Define a linear relation $T$ on $\ell^{2}$ by $$T(0, x_{2}, x_{3}, \dots,x_{n},\dots) = \{(z, 2{x_{2}}, 3{x_{3}},\dots, n{x_{n}},\dots) : z \in \mathbb{C}\},$$ with domain $$D(T) = \{(0, x_{2}, x_{3}, \dots,x_{n},\dots): (0, 2{x_{2}}, 3{x_{3}},\dots, n{x_{n}},\dots) \in \ell^{2}, x_{i}\in \mathbb{C}, i \in \mathbb{N}\setminus \{1\}\}.$$ The relation $T$ is closed, Since it is the inverse of a self-adjoint operator $S: \ell^{2} \to \ell^{2}$ defined by $$S(x_{1}, x_{2},\dots,x_{n},\dots) = (0, \frac{x_{2}}{2},\dots, \frac{x_{n}}{n},\dots).$$ Moreover, $R(T)$ is closed because $R(T) = D(S) = \ell^{2}.$\\
 Finally, for any $(y_{1}, y_{2},\dots, y_{n},\dots) \in \ell^{2}$, the Moore-Penrose inverse of $T$ is given by 
 \begin{align*}
 T^{\dagger}(y_{1},y_{2},\dots,y_{n},\dots) &= P_{N(T)^{\perp}}T^{-1}P_{R(T)} (y_{1},y_{2},\dots,y_{n},\dots)\\ 
 & = S(y_{1},y_{2},\dots,y_{n},\dots)\\ 
 & = (0,\frac{y_{2}}{2},\dots,\frac{y_{n}}{n},\dots).
 \end{align*}
 \end{example}
 \begin{example}
 Define a linear relation $T$ on the domain $\mathbb{R}^{3}$ by $$T(x_{1},x_{2}, x_{3}) = \{(x_{1}, x_{2},z) : \text{ for all}~ z \in \mathbb{R}\}.$$ The multivalued part of $T$ is $$M(T) =\{(0,0, z) : \text{for all}~z\in \mathbb{R}\}$$ and for each $(x_{1}, x_{2}, x_{3}) \in \mathbb{R}^{3},$ $$\|T(x_{1}, x_{2}, x_{3})\| = \|(x_{1},x_{2}, 0)\|.$$ Hence, $T$ is closed (in fact, bounded) and its range satisfies $R(T) = \mathbb{R}^{3}.$\\
 Moreover, $$(N(T))^{\perp} = \{(z_{1}, z_{2}, 0): z_{1}, z_{2} \in \mathbb{R}\},$$ from which it follows that the Moore–Penrose inverse of $T$ is given by $$T^{\dagger}(y_{1}, y_{2}, y_{3}) = (y_{1}, y_{2}, 0)~ \text{for all} ~ (y_{1}, y_{2}, y_{3}) \in \mathbb{R}^{3}.$$ In particular, $(T^{\dagger})^{*} = T^{\dagger}.$ Next, the adjoint relation $T^{*}$ is defined on $$D(T^{*}) =\{(z_{1},z_{2}, 0):(z_{1},z_{2}, 0) \in \mathbb{R}^{3}\},$$ and satisfies $$T^{*}(z_{1}, z_{2}, 0) = (z_{1}, z_{2},0).$$ Consequently, for every $(w_{1}, w_{2}, w_{3}) \in \mathbb{R}^{3},$ $$(T^{*})^{\dagger}(w_{1}, w_{2}, w_{3}) = (w_{1}, w_{2}, 0) = T^{\dagger}(w_{1}, w_{2}, w_{3}) = (T^{\dagger})^{*}(w_{1}, w_{2}, w_{3}).$$ This example therefore provides a concrete verification of Theorem \ref{thm 3.1}.
\end{example}
 \begin{theorem}\label{thm 3.3}
 Let $T \in CR(H, K)$. Then the following statements hold:
 \begin{enumerate}
 \item $R(T^{\dagger}) = (N(T))^{\perp} \cap D(T).$
 \item $N(T^{\dagger}) = (R(T))^{\perp} + M(T).$
 \item $TT^{\dagger}T = T.$
 \item $T^{\dagger}TT^{\dagger} = T^{\dagger}.$
 
  \item $(T^{\dagger})^{\dagger} \subset T$ on $D(T).$
  \end{enumerate}
  Moreover, if $R(T)$ is closed, then:
  \begin{enumerate}
  \setcounter{enumi}{5}
 \item $(T^{*}T)^{\dagger} = T^{\dagger}(T^{*})^{\dagger}.$
 \item $(TT^{*})^{\dagger} = (T^{*})^{\dagger}T^{\dagger}.$

 \end{enumerate}
 \end{theorem}
 \begin{proof}
 $\emph{(1)}$ Let $x \in R(T^{\dagger}).$ Then there exists $y \in K$ such that $\{y, x\} \in T^{\dagger},$ Writing $y = y^{'} + y^{''}$, with $y^{'} \in R(T)$, $y^{''} \in (R(T))^{\perp},$ we have $$\{y, y^{'}\} \in P_{R(T)},~~ \{y^{'}, x\} \in T^{-1},~~ \{x, x\} \in (N(T))^{\perp}.$$ Thus $x \in (N(T))^{\perp} \cap D(T)$ showing $$R(T^{\dagger}) \subset (N(T))^{\perp} \cap D(T).$$ Conversely, let $z \in (N(T))^{\perp} \cap D(T)$. Then there exists $v \in R(T)$ such that $$\{z, v\} \in T, ~~\{v, v\} \in P_{R(T)},~~ \{z, z\} \in P_{(N(T))^{\perp}},$$ which implies $\{v, z\} \in T^{\dagger}.$ Hence, $z \in R(T^{\dagger}).$ So, $$(N(T))^{\perp} \cap D(T) \subset R(T^{\dagger}).$$ Therefore, $(N(T))^{\perp} \cap D(T) = R(T^{\dagger}).$\\
 
 \medskip
\emph{(2)} Let $u \in (R(T))^{\perp}.$ Then $\{u, 0\} \in P_{(N(T))^{\perp}} T^{-1}P_{R(T)} = T^{\dagger},$ and hence  $u \in N(T^{\dagger}).$ Thus, $$(R(T))^{\perp} \subset N(T^{\dagger}).$$ If $k \in M(T) \subset R(T),$ then $\{k, 0\} \in T^{-1}$. So, $\{k, 0\} \in T^{\dagger}$. Hence, $$(R(T))^{\perp} + M(T) \subset N(T^{\dagger}).$$ Conversely, let $w \in N(T^{\dagger}),$ and decompose $$w = w^{'} + w^{''} ~\text{with}~ w^{'} \in R(T),~ w^{''} \in (R(T))^{\perp}.$$ Then, $$\{w, w^{'}\} \in P_{R(T)}, ~~ \{w^{'}, 0\} \in T^{-1},~~ \{0, 0\} \in P_{(N(T))^{\perp}},$$ which shows $w = w^{'} + w^{''} \in M(T) + (R(T))^{\perp}.$ Therefore, $N(T^{\dagger}) = (R(T))^{\perp} + M(T).$\\
 
 $\emph{(3)}$ Let $\{x, y\} \in TT^{\dagger}T$. Then there exist $z \in R(T)$ and $w \in D(T) \cap (N(T))^{\perp}$ such that $$\{x, z\} \in T, ~~\{z, w\} \in T^{\dagger}, ~~\{w, y\} \in T.$$ So, $\{z, z\} \in P_{R(T)}$, $\{z, w\} \in T^{-1}$, $\{w, w\} \in P_{(N(T))^{\perp}}.$ This yields $\{0, y-z\} \in T,$ and hence  $\{x, y\} = \{x, z\} + \{0, y-z\} \in T,$ proving $$TT^{\dagger}T \subset T.$$\\
 Let $\{u, v\} \in T$. Then $$\{v, v\} \in P_{R(T)},~ \{v, u\} \in T^{-1},~ \{u, u^{'}\} \in P_{(N(T))^{\perp}},$$ where $u = u^{'} + u^{''}$ and $u^{'} \in D(T) \cap (N(T))^{\perp}$, $u^{''} \in N(T).$\\ So, $\{v, u^{'}\} \in T^{\dagger}$ and $\{u^{'}, v\} \in T$ which shows that $\{u, v\} \in TT^{\dagger}T.$ Therefore, $$TT^{\dagger}T = T.$$
 
 $\emph{(4)}$ Let $\{x, y\} \in T^{\dagger}TT^{\dagger}$. Then there exist $z \in R(T^{\dagger}) = (N(T))^{\perp} \cap D(T)$ and $w \in R(T)$ such that $$\{x, z\} \in T^{\dagger},~~ \{z, w\} \in T,~~ \{w, y\} \in T^{\dagger}.$$ Writing $x = x^{'} + x^{''}$, with $x^{'} \in R(T)$, $x^{''} \in (R(T))^{\perp},$ we obtain $$\{x, x^{'}\} \in P_{R(T)},~~ \{x^{'}, z\} \in T^{-1},~~ \{z, z\} \in P_{(N(T))^{\perp}},$$ and $$\{w, w\} \in P_{R(T)},~~ \{w, y\} \in T^{-1},~~ \{y, y\} \in P_{(N(T))^{\perp}}.$$ Thus, $\{z-y, 0\} \in T$ which implies  $z- y \in N(T) \cap (N(T))^{\perp} = \{0\}.$ So, $\{x, y\} \in T^{\dagger}.$ Hence, $T^{\dagger}TT^{\dagger} \subset T^{\dagger}.$\\
  To prove the reverse inclusion, consider $\{u, v\} \in T^{\dagger}$ with $u = u^{'} + u^{''}$, and $u^{'} \in R(T)$, $u^{''} \in (R(T))^{\perp}.$ Thus, $$\{u, u^{'}\} \in P_{R(T)},~ \{u^{'}, v\} \in T^{-1},~ \{v, v\} \in P_{(N(T))^{\perp}},$$ and $$\{v, u^{'}\} \in T,~ \{u^{'}, u^{'}\} \in P_{R(T)},~ \{u^{'}, v \} \in T^{-1},~ \{v, v\} \in P_{(N(T))^{\perp}}.$$ This shows that $\{u, v\} \in T^{\dagger}TT^{\dagger}$. Hence, $T^{\dagger} \subset T^{\dagger}TT^{\dagger}$. Therefore, $T^{\dagger}TT^{\dagger} = T^{\dagger}.$\\
  
  $\emph{(5)}$ On $D(T)$, the Moore–Penrose inverse of $T^{\dagger}$ satisfies $$(T^{\dagger})^{\dagger} = P_{(N(T^{\dagger}))^{\perp}} (T^{\dagger})^{-1} (P_{R(T^{\dagger})}\vert_{D(T)}).$$ Let $\{x, y\} \in (T^{\dagger})^{\dagger}$ where $x \in D(T).$ Using the decompositions of $x = x^{'} + x^{''}$ with $x^{'} \in (N(T))^{\perp} \cap D(T)$, $x^{''} \in N(T).$ Then $$\{x, x^{'}\} \in P_{(N(T))^{\perp}},~ \{x^{'}, y\} \in (T^{\dagger})^{-1},~ \{y, y\} \in P_{(N(T^{\dagger}))^{\perp}}.$$
  The element $\{y, x^{'}\} \in T^{\dagger}$ and $y \in R(T)$ confirm that $$\{y, y\} \in P_{R(T)},~ \{y, x^{'}\} \in T^{-1},~ \{x^{'}, x^{'}\} \in P_{(N(T))^{\perp}}.$$ Hence, $\{x, y\} = \{x^{'}, y\} + \{x^{''}, 0\} \in T$. Therefore, $(T^{\dagger})^{\dagger} \subset T$ on $D(T)$.\\
  
  $\emph{(6)}$ Assume that $R(T)$ is closed. First, we will show that $(T^{*}T)^{\dagger} \subset T^{\dagger}(T^{*})^{\dagger}.$ By Proposition 2.4 \cite{MR3079830}, $$(T^{*}T)^{\dagger} = P_{(N(T))^{\perp}}(T^{*}T)^{-1}P_{R(T^{*})} = P_{(N(T))^{\perp}} T^{-1}(T^{*})^{-1}P_{R(T^{*})}.$$ Let $\{x, y\} \in (T^{*}T)^{\dagger}.$ Writing $x = x^{'} + x^{''}$ and $x^{'} \in R(T^{*})$, $x^{''} \in N(T).$ There exists $z \in D(T^{*}) \cap (N(T^{*}))^{\perp}$ such that $$\{x, x^{'}\} \in P_{R(T^{*})},~ \{x^{'}, z\} \in (T^{*})^{-1},~ \{z, y\} \in T^{-1},~\{y, y\} \in P_{(N(T))^{\perp}}.$$ So, $\{x, x^{'}\} \in P_{R(T^{*})}$, $\{x^{'}, z\} \in (T^{*})^{-1}$, $\{z, z\} \in P_{(N(T^{*}))^{\perp}}$, $\{z, z\} \in P_{R(T)}$, $\{z, y\} \in T^{-1}$, $\{y, y\} \in P_{(N(T))^{\perp}}.$\\ Thus, $\{x, y\} \in T^{\dagger} (T^{*})^{\dagger}.$ Hence,  $(T^{*}T)^{\dagger} \subset T^{\dagger}(T^{*})^{\dagger}.$\\
  To prove the reverse inclusion, consider $\{u, v\} \in T^{\dagger}(T^{*})^{\dagger}$. Then there exists $w \in R((T^{*})^{\dagger}) = D(T^{*}) \cap (N(T^{*}))^{\perp}$ such that $\{u, w\} \in (T^{*})^{\dagger}$ and $\{w, v\} \in T^{\dagger}.$\\ Writing $u = u^{'} + u^{''}$ and $u^{'} \in R(T^{*})$, $u^{''} \in N(T),$ we get, $$\{u, u^{'}\} \in P_{R(T^{*})},~ \{u^{'}, w\} \in (T^{*})^{-1},~ \{w, w\} \in P_{R(T)},~ \{w, v\} \in T^{-1},~ \{v, v\} \in P_{(N(T))^{\perp}}.$$ So, $\{u, v\} \in (T^{*}T)^{\dagger}$. Hence, $T^{\dagger}(T^{*})^{\dagger} \subset (T^{*}T)^{\dagger}.$ Therefore, $T^{\dagger}(T^{*})^{\dagger} = (T^{*}T)^{\dagger}.$\\
  
  $\emph{(7)}$ To prove $(TT^{*})^{\dagger} = (T^{*})^{\dagger}T^{\dagger}$, we replace $T$ by $T^{*}$ in part~(6) yields.
  
\end{proof}
Proposition 3.1 \cite{MR3079830} says that the reduced minimal modulus of $T$, $\gamma(T) = \frac{1}{\|T^{\dagger}\|}$ when $T \in CR(H, K)$ with closed range and $T^{\dagger} = (T^{-1})_{op}.$ Theorem \ref{thm 3.4} depicts the generalized version of Proposition 3.1 \cite{MR3079830}.
\begin{theorem}\label{thm 3.4}
Let $T \in CR(H, K)$ be a closed range linear relation. Then $$\gamma(T) = \frac{1}{\|T^{\dagger}\|}.$$
\end{theorem}
\begin{proof}
By Theorem \ref{thm 3.2}, the Moore–Penrose inverse $T^{\dagger}$ is an operator; in particular, $$\|T^{\dagger}\| < \infty.$$
\noindent
Case $(1):$ Consider $N(T^{\dagger}) = (R(T))^{\perp} + M(T) \neq K$. Then,
\begin{align*}
\|T^{\dagger}\| &= \sup\biggl\{\frac{\|x\|}{\|y\|} : \{y, x\} \in T^{\dagger} ~(\text{Here,}~ T^{\dagger}{y} = x), x \neq 0\biggr\}\\
&=\frac{1}{\inf\biggl\{\frac{\|y\|}{\|x\|} : \{y, x\} \in T^{\dagger}, x \neq 0\biggr\}}\\
&=\frac{1}{\inf\biggl\{\frac{\|y^{'}\|}{\|x\|} : \{y, x\} \in T^{\dagger}, x \neq 0, \{y,y^{'}\} \in P_{R(T)}\biggr\}}\\
&=\frac{1}{\inf\biggl\{\frac{\|y^{'}\|}{\|x\|} : \{x, y^{'}\} \in T, x \in (N(T))^{\perp}\setminus\{0\} \biggr\}}\\
&= \frac{1}{\gamma(T)} ~(\text{By Proposition II.2.2 \cite{MR1631548}}).
\end{align*}
Case $(2):$ When $N(T^{\dagger}) = R(T)^{\perp} + M(T) = K$, then $\|T^{\dagger}\| = 0$. Let $h \in D(T).$ Then there exists $k \in K$ such that $\{h, k\} \in T.$ Decompose $k = k^{'} + k^{''}$, where $k^{'} \in R(T)^{\perp}$, $k^{''} \in M(T).$ So, $\{h, k^{'}\} \in T$ and $$k^{'} \in R(T) \cap (R(T))^{\perp} = \{0\},$$ it follows that $h \in N(T)$ and $D(T) \subset N(T).$ By Proposition II.2.2 \cite{MR1631548}, we have, $$\gamma(T) = \infty = \frac{1}{\|T^{\dagger}\|}.$$
\end{proof}
Proposition 4.1(ii) \cite{MR3079830} says that $\lambda \in \sigma(T) \setminus \{0\} \text{ if and only if } \frac{1}{\lambda} \in \sigma(T^{\dagger}) \setminus \{0\},$ when $T \in CR(H, K)$ with closed range and $T^{\dagger} = (T^{-1})_{op}.$ Theorem \ref{thm 3.5} depicts the generalized version of Proposition 4.1(ii) \cite{MR3079830}.
\begin{theorem}\label{thm 3.5}
Let $T \in CR(H, H)$ be self-adjoint. Then
\begin{align*}
\lambda \in \sigma(T) \setminus \{0\} \text{ if and only if } \frac{1}{\lambda} \in \sigma(T^{\dagger}) \setminus \{0\}. 
\end{align*}
\end{theorem}
\begin{proof}
\noindent
Assume that $\lambda \in \sigma(T) \setminus \{0\}$ and suppose, to the contrary, that $\frac{1}{\lambda} \in \rho(T^{\dagger})$. Then the bounded operator $(T^{\dagger} - \frac{1}{\lambda})$ is bijective, and defined on all of $H$.\\ Let $\{0, k\} \in (T - \lambda)^{-1}.$ Then $\{k, \lambda k\} \in T,$ which implies $k \in R(T) \cap D(T).$ Moreover, $$\{\lambda k, \lambda k\} \in P_{R(T)},~ \{\lambda k, k\} \in T^{-1},~ \{k, k\} \in P_{(N(T))^{\perp}}.$$ So, $\{\lambda k, k\} \in T^{\dagger}$ which implies  $$\{\lambda k, 0\} \in (T^{\dagger} - \frac{1}{\lambda}).$$ We conclude that $k = 0.$ Thus, $(T - \lambda)^{-1}$ is an operator. Since $(T - \lambda)$ is closed, its inverse $(T - \lambda)^{-1}$ is also closed. We now show that $D((T - \lambda)^{-1}) = H.$\\
Let $h \in H.$ Then there exists $v \in H$ such that $$(T^{\dagger} - \frac{1}{\lambda})v = h,$$ or equivalently, $$T^{\dagger}v = \frac{1}{\lambda}v + h \in R(T^{\dagger}) = (N(T))^{\perp} \cap D(T).$$ Write $v = v^{'} + v^{''}$, where $v^{'} \in R(T)$ and $v^{''} \in (R(T))^{\perp}.$ Then $$\{v, v^{'}\} \in P_{R(T)},~ \{v^{'}, \frac{1}{\lambda} v + h\} \in T^{-1},~ \{\frac{1}{\lambda} v + h, \frac{1}{\lambda} v + h\} \in P_{(N(T))^{\perp}}.$$ Thus, $\{\frac{1}{\lambda}v + h, v^{'}\} \in T$ together with $\{v^{''}, 0\} \in T,$ we obtain $$\{\frac{1}{\lambda}v^{'} + h, - \lambda h\} \in T- \lambda.$$ Hence, $h \in D((T - \lambda)^{-1}).$ By Theorem III.4.2 \cite{MR1631548}, we get $(T - \lambda)^{-1}$ is bounded operator in domain $H$ which is  a contradiction. Our assumption is wrong. Therefore, $$\frac{1}{\lambda} \in \sigma(T^{\dagger}) \setminus \{0\}.$$

To prove converse part, we assume that $\frac{1}{\mu} \in \sigma(T^{\dagger}) \setminus \{0\},$ and suppose that $\mu \in \rho(T).$ Then $(T - \mu)^{-1}$ is a bounded operator on domain $H$. Let $x \in N(T^{\dagger} - \frac{1}{\mu}).$ Then $$T^{\dagger}x = \frac{1}{\mu} x \in (N(T))^{\perp} \cap D(T) \cap D(T^{\dagger}) = R(T) \cap D(T).$$ Thus, $$\{x, x\} \in P_{R(T)},~ \{x, \frac{1}{\mu}x\} \in T^{-1},~ \{\frac{1}{\mu}x, \frac{1}{\mu}x\} \in P_{(N(T))^{\perp}}.$$ Hence, $\{\frac{1}{\mu}x, x\} \in T,$ therefore $$\{\frac{1}{\mu}x, 0\} \in T - \mu.$$ Since $(T- \mu)^{-1}$ is an operator, it follows that $x \in N(T - \mu) = \{0\},$ and thus $(T^{\dagger} - \frac{1}{\mu})$ is injective. To prove that $(T^{\dagger} - \frac{1}{\mu})$ surjective, let $w \in H = D((T - \mu)^{-1}).$ Then there exists $s \in H$ such that $\{s, w\} \in (T - \mu).$ Write $s = s^{'} + s^{''},$ where $s^{'} \in (N(T))^{\perp} \cap D(T)$, $s^{''} \in N(T).$ Then $\{s^{'}, w + \mu s\} \in T.$  Consequently, $$\{w + \mu s, w + \mu s\} \in P_{R(T)},~ \{w + \mu s, s^{'}\} \in T^{-1},~ \{s^{'}, s^{'}\} \in P_{(N(T))^{\perp}}.$$ Hence, $\{w + \mu s, s^{'}\} \in T^{\dagger}$ and $\{s^{''}, 0\} \in T^{\dagger}$ confirm that $\{w + \mu s^{'}, -\frac{1}{\mu}w\} \in (T^{\dagger} - \frac{1}{\mu}).$ Thus, $w \in R(T^{\dagger} - \frac{1}{\mu})$ and $(T^{\dagger} - \frac{1}{\mu})$ is surjective.. Therefore, $\frac{1}{\mu}\in \rho(T^{\dagger})$ is again a contradiction. Thus, $$\mu \in \sigma(T) \setminus \{0\}.$$
\end{proof}
\begin{theorem}\label{thm 3.6}
Let $T \in CR(H, K) \cap BR(H, K).$ Then $T^{\dagger}T$ is a closed operator. Moreover, if $R(T)$ is closed, then $TT^{\dagger}$ is a bounded linear relation.
\end{theorem}
\begin{proof}
We first show that $T^{\dagger}T$ is an operator. Suppose $\{0,k\} \in T^{\dagger}T.$ Then there exists $p \in R(T)$ such that $$\{0, p\} \in T ~~ \{p, k\} \in T^{\dagger}.$$ Consequently, $$\{p, p\} \in P_{R(T)},~ \{p, k\} \in T^{-1},~ \{k, k\} \in P_{(N(T))^{\perp}}.$$ It follows that $$\{k, 0\} = \{k, p\} - \{0, p\} \in T,$$ which implies $k \in N(T) \cap (N(T))^{\perp} = \{0\}.$ Hence, $T^{\dagger}T$ is an operator.\\
Let $\{x, y\} \in \overline{T^{\dagger}T}.$ Then there exists a sequence $\{\{x_{n}, y_{n}\}\}$ in $T^{\dagger}T$ such that $$\{x_{n}, y_{n}\} \to \{x, y\},~~ \text{as}~n \to \infty.$$ There is a sequence $\{z_{n}\}$ in $R(T)$ such that $\{x_{n}, z_{n}\} \in T$ and $\{z_{n}, y_{n}\} \in T^{\dagger}$, for all $n \in \mathbb{N}.$ Thus, $$\{z_{n}, z_{n}\} \in P_{R(T)},~ \{z_{n}, y_{n}\} \in T^{-1},~ \{y_{n}, y_{n}\} \in P_{(N(T))^{\perp}},$$ which yields $y_{n} \in D(T) \cap (N(T))^{\perp}$, for all $n \in \mathbb{N}.$ Since $T$ is continuous, the set $(N(T))^{\perp} \cap D(T)$ is closed. So, $$y \in (N(T))^{\perp} \cap D(T).$$ Then there exists $w \in R(T)$ such that $\{y, w\} \in T$. Moreover, $$\{w, w\} \in P_{R(T)},~ \{w, y\} \in T^{-1},~ \{y, y\} \in P_{(N(T))^{\perp}},$$ which yields $\{w, y\} \in T^{\dagger}.$ Since $\{x_{n}, z_{n}\} \in T$ and $\{y_{n}, z_{n}\} \in T,$ passing to the limit gives  $\{x-y, 0\} \in T.$ Hence, $\{x, w\} = \{x-y, 0\} + \{y, w\} \in T.$ Therefore, $\{x, y\} \in T^{\dagger}T.$ This proves that $T^{\dagger}T$ is a closed operator.\\ 
By Theorem \ref{thm 3.2}, $T^{\dagger}$ is a bounded operator. The composition $TT^{\dagger}$ is a closed linear relation. By Theorem III.4.2 \cite{MR1631548}, we obtain $TT^{\dagger}$ is continous in domian $K,$ and hence $TT^{\dagger}$ is bounded.

\end{proof}
\begin{theorem}\label{thm 3.7}
Let $T \in CR(H, K)$ and $S \in LR(H, K)$ satisfy $M(S) \subset M(T).$ Then the product $T^{\dagger}S$ is an operator. 
\end{theorem}
\begin{proof}
Let $\{0, k\} \in T^{\dagger}S$. Then there exists $w \in R(S)$ with $\{0, w\} \in S$ and $\{w, k\} \in T^{\dagger}$ such that $$\{w, w\} \in P_{R(T)}~ (\text{ because }  w \in M(S) \subset M(T)), ~\{w, k\} \in T^{-1},~ \{k, k\} \in P_{(N(T))^{\perp}}.$$ It follows that $$\{k, 0\} = \{k, w\} - \{0, w\} \in T,$$ which implies $$k \in N(T) \cap N(T)^{\perp} = \{0\}.$$ Therefore, $T^{\dagger}S$ is an operator.
\end{proof}
\begin{theorem}\label{thm 3.8}
Let $S \in CR(H, K)$ and $T \in CR(K, H_{1}).$ Suppose that $$R(S) = (N(T))^{\perp}.$$ Then $$(TS)^{\dagger} = S^{\dagger}T^{\dagger}.$$
\end{theorem}
\begin{proof}
We begin by identifying the range and kernel of the product $TS.$
\medskip
\noindent
 It is clear that $R(TS) \subset R(T).$ Conversely, let $y \in R(T)$. Then there exists $x \in (N(T))^{\perp} \cap D(T)$ such that $\{x, y\} \in T.$ Since $R(S) = (N(T))^{\perp},$ there exists $w \in D(S)$ such that $\{w, x\} \in S.$ Thus, $\{w, y\} \in TS$ which implies $y \in R(TS)$. Hence, $$R(TS) = R(T).$$ 
 Clearly $N(S) \subset N(TS).$ To prove the reverse inclusion, let $z \in N(TS).$ Then there exists $p \in R(S)$ such that $$\{z, p\} \in S, ~ \{p, 0\} \in T.$$ Thus, $$p \in N(T) \cap (N(T))^{\perp} = \{0\},$$ which implies that $z \in N(S).$ Hence, $$N(TS) = N(S).$$ Using the characterizations of the Moore–Penrose inverse, we write $$(TS)^{\dagger} = P_{(N(S))^{\perp}} S^{-1}T^{-1} P_{R(T)}.$$ Let $\{u, v\} \in (TS)^{\dagger}$, and decompose $u = u^{'} + u^{''}$ with $u^{'} \in R(T)$, $u^{''} \in (R(T))^{\perp}.$ Then there exist $r \in D(T) \cap (N(T))^{\perp}$ and $v \in D(S)$ such that $$\{u, u^{'}\} \in P_{R(T)},~ \{u^{'}, r\} \in T^{-1},~ \{r, v\} \in S^{-1},~ \{v, v\} \in P_{(N(S))^{\perp}}.$$ Since $R(S) = (N(T))^{\perp},$ we have $$\{u, u^{'}\} \in P_{R(T)},~ \{u^{'}, r\} \in T^{-1},~ \{r, r\} \in P_{(N(T))^{\perp}},$$ and $$\{r, r\} \in P_{R(S)},~ \{r, v\} \in S^{-1}, \{v, v\} \in P_{(N(S))^{\perp}}.$$ Consequently, $\{u, v\} \in S^{\dagger} T^{\dagger},$ and hence $(TS)^{\dagger} \subset S^{\dagger}T^{\dagger}.$\\
To prove the reverse inclusion, let $\{h, k\} \in S^{\dagger}T^{\dagger} = P_{(N(S))^{\perp}} S^{-1} P_{R(S)}T^{-1}P_{R(T)}.$ Then there exists $q \in D(T) \cap (N(T))^{\perp}$  and $h = h^{'} + h^{''}$  with $h^{'} \in R(T)$, $h^{''} \in (R(T))^{\perp}$ such that $$\{h, h^{'}\} \in P_{R(T)},~ \{h^{'}, q\} \in T^{-1},~ \{q, q\} \in P_{(N(T))^{\perp}} = P_{R(S)},~ \{q, k\} \in S^{-1},~ \{k, k\} \in P_{(N(S))^{\perp}}.$$ These relations imply that $\{h, k\} \in (TS)^{\dagger},$ hence  $$S^{\dagger}T^{\dagger} \subset (TS)^{\dagger}.$$ Combining the two inclusions yields $$(TS)^{\dagger} = S^{\dagger} T^{\dagger}.$$
\end{proof}
\begin{theorem}\label{thm 3.9}
Let $T \in CR(H, K)$ be a closed range linear relation. Then $$(T^{*}T)^{\dagger}T^{*} \subset T^{\dagger} \subset T^{*}(TT^{*})^{\dagger}.$$ Moreover, the inclusion on the right becomes an equality whenever $D(T)$ is dense in $H.$
\end{theorem}
\begin{proof}
By Theorem \ref{thm 3.3}, it suffices to establish the inclusions $$T^{\dagger}(T^{*})^{\dagger}T^{*} \subset T^{\dagger} \subset T^{*}(T^{*})^{\dagger} T^{\dagger}.$$ First, let $\{x, y\} \in T^{\dagger}(T^{*})^{\dagger}T^{*}.$ Using the representation $$T^{\dagger}(T^{*})^{\dagger}T^{*} = P_{(N(T))^{\perp}} T^{-1}P_{R(T)} (T^{*})^{-1} P_{R(T^{*})} T^{*} = P_{(N(T))^{\perp}} T^{-1}P_{R(T)} (T^{*})^{-1} T^{*},$$ there exist elements $z \in R(T^{*})$, $w \in D(T^{*})$ such that $w = w^{'} + w^{''}$, where $w^{'} \in R(T) \cap D(T^{*})$, $w^{''} \in (R(T))^{\perp}$ satisfying $$\{x, z\} \in T^{*},~ \{z, w\} \in (T^{*})^{-1},~ \{w, w^{'}\} \in P_{R(T)},~ \{w^{'}, y\} \in T^{-1},~ \{y, y\} \in P_{(N(T))^{\perp}}.$$ 
 Writing $x = x^{'} + x^{''}$, with $x^{'} \in (N(T^{*}))^{\perp} \cap D(T^{*})$ and $x^{''} \in N(T^{*}),$ we obtain $$\{x^{'}, z\} \in T^{*}, ~ \{w^{'}, z\} \in T^{*},$$ which implies $$x^{'}- w^{'} \in N(T^{*}) \cap (N(T^{*}))^{\perp} = \{0\}.$$ Hence, $$\{x, x^{'}\} \in P_{R(T)},~ \{x^{'}, y\} \in T^{-1},~ \{y, y\} \in P_{(N(T))^{\perp}},$$ and therefore $\{x, y\} \in T^{\dagger},$ showing that
\begin{align}\label{equ 1}
T^{\dagger}(T^{*})^{\dagger}T^{*} \subset T^{\dagger}.
\end{align}
Next, let $\{u, v\} \in T^{\dagger} = P_{(N(T))^{\perp}}T^{-1}P_{R(T)},$ where $u = u^{'} + u^{''}$ with $u^{'} \in R(T)$ and $u^{''} \in (R(T))^{\perp}.$ Then $$\{u, u^{'}\} \in P_{R(T)},~ \{u^{'}, v\} \in T^{-1},~ \{v, v\} \in P_{(N(T))^{\perp}}.$$ Since $v \in R(T^{*}),$ there exists $s \in D(T^{*}) \cap (N(T^{*}))^{\perp}$ such that $\{s, v\} \in T^{*}$. Thus, $$\{u, u^{'}\} \in P_{R(T)},~ \{u^{'}, v\} \in T^{-1},~ \{v, v\} \in P_{R(T^{*})},~ \{v, s\} \in (T^{*})^{-1},$$ $$\{s, s\} \in P_{(N(T^{*}))^{\perp}},~ \{s, v\} \in T^{*}.$$ This shows that $\{u,v\} \in T^{*}(T^{*})^{\dagger}T^{\dagger}.$ Hence,
\begin{align}\label{equ 2}
T^{\dagger} \subset T^{*}(T^{*})^{\dagger}T^{\dagger}.
\end{align}
Combining the inclusions (\ref{equ 1}) and (\ref{equ 2}), we obtain 
\begin{align*}
(T^{*}T)^{\dagger}T^{*} \subset T^{\dagger} \subset T^{*}(TT^{*})^{\dagger}.
\end{align*}
Finally, assume that $D(T)$ is dense in $H.$ Let $$\{p, q\} \in T^{*}(T^{*})^{\dagger} T^{\dagger} = T^{*}P_{R(T)}(T^{*})^{-1}P_{R(T^{*})}T^{-1}P_{R(T)},$$ and write $p = p^{'} + p^{''}$ with $p^{'} \in R(T)$, $p^{''} \in (R(T))^{\perp}.$\\ Then there exist $r \in D(T) \cap (N(T))^{\perp}$ and $t \in D(T^{*}) \cap (N(T^{*}))^{\perp}$ such that $$\{p, p^{'}\} \in P_{R(T)},~ \{p^{'}, r\} \in T^{-1},~ \{r, r\} \in P_{R(T^{*})},~ \{r, t\} \in (T^{*})^{-1},~ \{t, t\} \in P_{R(T)},~ \{t, q\} \in T^{*}.$$ This yields $\{r- q, 0\} \in (T^{*})^{-1}$ which implies $r- q \in M(T^{*}) = (D(T))^{\perp} = \{0\}.$ Hence, $$\{p, p^{'}\} \in P_{R(T)},~ \{p^{'}, q\} \in T^{-1},~ \{q, q\} \in P_{(N(T))^{\perp}}$$ which shows $\{p, q\} \in T^{\dagger}.$ Therefore, $$T^{\dagger} = T^{*}(TT^{*})^{\dagger},$$ when $D(T)$ is dense in $H$.
\end{proof}
\begin{theorem}\label{thm 3.10}
Let $T \in CR(H, K)$ be a closed range linear relation. Then $$(T^{*})^{\dagger}T^{*}T \subset T \subset TT^{*}(T^{*})^{\dagger}.$$
\end{theorem}
\begin{proof}
Let $\{x, y\} \in  (T^{*})^{\dagger}T^{*}T.$ Using the representation $$(T^{*})^{\dagger}T^{*}T = P_{R(T)}(T^{*})^{-1}P_{R(T^{*})}T^{*}T,$$ there exist $z \in R(T)$ and $w \in R(T^{*})$ such that $$\{x, z\} \in T,~ \{z, w\} \in T^{*},~ \{w, w\} \in P_{R(T^{*})},~ \{w, y\} \in (T^{*})^{-1},~ \{y, y\} \in P_{R(T)}.$$ Then $\{z-y, 0\} \in T^{*},$ it follows that $$z-y \in N(T^{*}) \cap (N(T^{*}))^{\perp} = \{0\}.$$ Consequently $\{x, y\} = \{x, z\} \in T$. This proves that $$(T^{*})^{\dagger}T^{*}T \subset T.$$
Next, we prove the right-hand inclusion. Let $\{u, v\} \in T.$ Decompose $u = u^{'} + u^{''}$, where $u^{'} \in D(T) \cap (N(T))^{\perp}$ and $u^{''} \in N(T).$\\ There exists $s \in (N(T^{*}))^{\perp} \cap D(T^{*})$ with $\{s, u^{'}\} \in T^{*}$ such that $$\{u, u^{'}\} \in P_{R(T^{*})},~ \{u^{'}, s\} \in (T^{*})^{-1},~ \{s, s\} \in P_{R(T)},~ \{s, u^{'}\} \in T^{*},~ \{u^{'}, v\} \in T.$$ These relations show that $\{u, v\} \in TT^{*}(T^{*})^{\dagger}$. Therefore, $T \subset TT^{*}(T^{*})^{\dagger}.$
\end{proof}
\begin{theorem}\label{thm 3.11}
Let $T \in CR(H, K)$ be a closed range linear relation. Then $$(T^{\dagger})^{*}T^{\dagger}T \subset (T^{\dagger})^{*} \subset TT^{\dagger}(T^{\dagger})^{*}.$$ Moreover, the right-hand inclusion becomes an equality if $D(T^{*})$ is dense in $K.$
\end{theorem}
\begin{proof}
We begin with the left-hand inclusion. Let $\{x, y\} \in (T^{\dagger})^{*}T^{\dagger}T.$ Using the representation $$(T^{\dagger})^{*}T^{\dagger}T = P_{R(T)}(T^{*})^{-1}P_{R(T^{*})} T^{-1}P_{R(T)}T,$$ there exist $z \in R(T)$ and $w \in (N(T))^{\perp} \cap D(T)$ such that $$\{x, z\} \in T,~ \{z, z\} \in P_{R(T)},~ \{z, w\} \in T^{-1},~ \{w, w\} \in P_{R(T^{*})},~ \{w, y\} \in (T^{*})^{-1},~ \{y, y\} \in P_{R(T)}.$$ Write $x = x^{'} + x^{''}$, where $x^{'} \in (N(T))^{\perp} \cap D(T)$ and $x^{''} \in N(T).$ Then $$\{x^{'}, z\} \in T,~  \{z, z\} \in P_{R(T)},~ \{z, w\} \in T^{-1},~ \{w, w\} \in P_{R(T^{*})},~ \{w, y\} \in (T^{*})^{-1},~ \{y, y\} \in P_{R(T)}.$$ Again, $\{x^{'} - w, 0\} \in T$ which implies $x^{'}- w \in N(T) \cap (N(T))^{\perp} = \{0\}.$ This shows that $$\{x, x^{'}\} \in P_{R(T^{*})},~ \{x^{'}, y\} \in (T^{*})^{-1},~ \{y, y\} \in P_{(N(T^{*}))^{\perp}}.$$ Hence, $\{x, y\} \in (T^{*})^{\dagger} = (T^{\dagger})^{*},$ proving that $$(T^{\dagger})^{*}T^{\dagger}T \subset (T^{\dagger})^{*}.$$
Next, we prove the right-hand inclusion. Let $\{u, v\} \in (T^{*})^{\dagger} = P_{R(T)}(T^{*})^{-1} P_{R(T^{*})}$ with $u = u^{'} + u^{''}$ where $u^{'} \in R(T^{*})$, $u^{''} \in (R(T^{*}))^{\perp}.$ Again, $$\{u, u^{'}\} \in P_{R(T^{*})},~ \{u^{'}, v\} \in (T^{*})^{-1},~ \{v, v\} \in P_{R(T)}.$$ Since $v \in R(T)$, there exists $s \in (N(T))^{\perp} \cap D(T)$ such that $\{s, v\} \in T.$ Consequently, $$\{u, u^{'}\} \in P_{R(T^{*})},~ \{u^{'}, v\} \in (T^{*})^{-1},~ \{v, v\} \in P_{R(T)},~ \{v, s\} \in T^{-1},~ \{s, s\} \in P_{(N(T))^{\perp}},~ \{s, v\} \in T.$$ These relations show that $\{u, v\} \in TT^{\dagger}(T^{*})^{\dagger},$ hence $$(T^{\dagger})^{*} \subset TT^{\dagger}(T^{\dagger})^{*}.$$
Finally, assume that $D(T^{*})$ is dense in $K.$
Let $$\{p, q\} \in TT^{\dagger}(T^{\dagger})^{*} = TP_{(N(T))^{\perp}}T^{-1} P_{R(T)} (T^{*})^{-1}P_{R(T^{*})}.$$ Write $p = p^{'} + p^{''}$, where $p^{'} \in R(T^{*})$ and $p^{''} \in (R(T^{*}))^{\perp}$, we obtain $$r \in (N(T^{*}))^{\perp} \cap D(T^{*})~ \text{and}~ h \in (N(T))^{\perp} \cap D(T)$$ such that $$\{p, p^{'}\} \in P_{R(T^{*})},~ \{p^{'}, r\} \in (T^{*})^{-1},~ \{r, r\} \in P_{R(T)},~ \{r, h\} \in T^{-1},~ \{h, h\} \in P_{(N(T))^{\perp}},~ \{h, q\} \in T.$$
Now, it follows that $r-q \in M(T) = (D(T^{*}))^{\perp} = \{0\},$ which implies that $$\{p, p^{'}\} \in P_{R(T^{*})},~ \{p^{'}, q\} \in (T^{*})^{-1},~ \{q, q\} \in P_{R(T)} = P_{(N(T^{*}))^{\perp}},$$ and hence $\{p, q\} \in (T^{*})^{\dagger}.$ Therefore, $$(T^{\dagger})^{*} = TT^{\dagger}(T^{\dagger})^{*},$$ when $D(T^{*})$ is dense in $K$.
\end{proof}
\begin{theorem}\label{thm 3.12}
Let $T \in CR(H, K)$ be a closed range linear relation. Then $$T^{\dagger}TT^{*} \subset T^{*} = (T^{\dagger}T)^{*}T^{*}.$$
\end{theorem}
\begin{proof}
The left-hand inclusion follows immediately by applying Theorem \ref{thm 3.10} with $T^{*}$ in place of $T.$ By using Proposition 1.3.9 \cite{MR3971207}, we have
 $$(T^{\dagger}T)^{*}T^{*} = T^{*}(T^{*})^{\dagger}T^{*}$$ because $T^{\dagger}$ is a bounded operator and $T$ is closed.\\
 Let $\{x, y\} \in T^{*}(T^{*})^{\dagger}T^{*}.$ Then there exist $z \in R(T^{*})$ and $u \in D(T^{*}) \cap (N(T^{*}))^{\perp}$ with $x = x^{'} + x^{''}$, where $x^{'} \in (N(T^{*}))^{\perp} \cap D(T^{*})$, $x^{''} \in N(T^{*}),$ such that $$\{x^{'}, z\} \in T^{*},~ \{z, z\} \in P_{R(T^{*})},~ \{z, u\} \in (T^{*})^{-1},~ \{u, u\} \in P_{R(T)},~ \{u, y\} \in T^{*}.$$ So, $\{x^{'}- u, 0\} \in T^{*}$. Thus, $x^{'}-u \in N(T^{*}) \cap N(T^{*})^{\perp} = \{0\}$ which implies $$\{x, y\} = \{u, y\} + \{x^{''}, 0\} \in T^{*}.$$ Hence, $$(T^{\dagger}T)^{*}T^{*} \subset T^{*}.$$
 
To prove the reverse inclusion, let $\{u, v\} \in T^{*}.$ Then there exist $$r \in D(T^{*}) \cap (N(T^{*}))^{\perp}$$ such that $$\{u, v\} \in T^{*},~ \{v, v\} \in P_{R(T^{*})},~ \{v, r\} \in (T^{*})^{-1},~ \{r, r\} \in P_{(N(T^{*}))^{\perp}},~ \{r, v\} \in T^{*}.$$ Hence, $\{u, v\} \in T^{*}(T^{\dagger})^{*}T^{*}$. Therefore, $$T^{*} = (T^{\dagger}T)^{*}T^{*}.$$
\end{proof}

Let $H$ and $K$ be two Hilbert spaces. The space $H \bigoplus K$ defined by $$H \bigoplus K = \{(h, k): h\in H, k\in K\}$$ is a linear space with respect to addition and scalar multiplication defined by 
\begin{align*}
&~~~~~~~~~~~~~~~~~~~~~~(h_{1}, k_{1}) + (h_{2}, k_{2}) = (h_{1} + h_{2}, k_{1} + k_{2}), \text{ and }\\
&\lambda (h, k ) = (\lambda h, \lambda k), \text{ for all } h, h_{1}, h_{2}\in H, \text{ for all } k, k_{1}, k_{2}\in K \text{ and } \lambda\in \mathbb{K}, ~(\mathbb{K} = \mathbb{R} \text{ or } \mathbb{C}).
\end{align*}
Now, $H \bigoplus K$ is an inner product space with respect to the inner product given by
\begin{align*}
\langle (h_{1}, k_{1}) (h_{2}, k_{2}) \rangle = \langle h_{1}, h_{2}\rangle + \langle k_{1}, k_{2}\rangle, \text{ for all } h_{1}, h_{2}\in H, \text{ and } \text{ for all }k_{1}, k_{2} \in K. 
\end{align*}
The norm on $H \bigoplus K$ is defined by 
\begin{align*}
\|(h, k)\| = (\|h\|^{2} + \|k\|^{2})^{\frac{1}{2}}, \text{ for all } (h,k) \in H \times K.
\end{align*}
Moreover, The direct sum of two linear relations $T_{1} \text{ and } T_{2}$ from $D(T_{1}) \subset H_{1}$ to $K_{1}$ and from $D(T_{2}) \subset H_{2}$ to $K_{2}$ respectively is defined by 
\begin{align*}
T_{1} \bigoplus T_{2} = \bigl\{\{(x_{1}, x_{2}), (y_{1}, y_{2})\} : \{x_{1}, y_{1}\} \in T_{1}, \{x_{2}, y_{2}\} \in T_{2}\bigr\}
\end{align*}
from $D(T_{1}) \bigoplus D(T_{2}) \subset H_{1} \bigoplus H_{2}$ into $K_{1} \bigoplus K_{2}.$ This is same as considering the matrix linear relation $\begin{bmatrix} 
T_{1} & 0\\
0 & T_{2} \end{bmatrix}.$
\begin{theorem}\label{thm 3.15}
Let $T_{1} \in CR(H_{1}, K_{1})$ and $T_{2} \in CR(H_{2}, K_{2}).$ Then $$(T_{1} \bigoplus T_{2})^{\dagger} = T_{1}^{\dagger} \bigoplus T_{2}^{\dagger}.$$
\end{theorem}
\begin{proof}
It is straightforward to verify that $$R(T_{1} \bigoplus T_{2}) = R(T_{1}) \bigoplus R(T_{2}), ~~ N(T_{1} \bigoplus T_{2}) = N(T_{1}) \bigoplus N(T_{2}),$$ and consequently, $$(N(T_{1} \bigoplus T_{2}))^{\perp} = (N(T_{1}))^{\perp} \bigoplus (N(T_{2}))^{\perp}.$$ Therefore, $$(T_{1} \bigoplus T_{2})^{\dagger} = P_{(N(T_{1}))^{\perp} \bigoplus (N(T_{2}))^{\perp}} (T_{1} \bigoplus T_{2})^{-1} P_{R(T_{1}) \bigoplus R(T_{2})},$$ while $$T_{1}^{\dagger} \bigoplus T_{2}^{\dagger} = P_{(N(T_{1}))^{\perp}} T_{1}^{-1} P_{R(T_{1})} \bigoplus P_{(N(T_{2}))^{\perp}} T_{2}^{-1} P_{R(T_{2})}.$$
Let $\{x, y\} \in T_{1}^{\dagger} \bigoplus T_{2}^{\dagger}.$ Then there exists $\{x_{1}, y_{1}\} \in T_{1}^{\dagger}$ and $\{x_{2}, y_{2}\} \in T_{2}^{\dagger}$ such that $\{x, y\} = \{(x_{1}, x_{2}), (y_{1}, y_{2})\}.$  Decompose $$x_{1} = x_{1}^{'} + x_{1}^{''}, ~\text{and}~ x_{2} = x_{2}^{'} + x_{2}^{''},$$ where $x_{1}^{'} \in R(T_{1}), x_{1}^{''} \in (R(T_{1}))^{\perp} , x_{2}^{'} \in R(T_{2}), \text{ and } x_{2}^{''} \in (R(T_{2}))^{\perp}.$ Then $$\{x_{1}, x_{1}^{'}\} \in P_{R(T_{1})},~ \{x_{1}^{'}, y_{1}\} \in T_{1}^{-1},~ \{y_{1}, y_{1}\} \in P_{(N(T_{1}))^{\perp}},$$ and $$\{x_{2}, x_{2}^{'}\} \in P_{R(T_{2})},~ \{x_{2}^{'}, y_{2}\} \in T_{2}^{-1},~ \{y_{2}, y_{2}\} \in P_{(N(T_{2}))^{\perp}}.$$
So, $\{(x_{1}, x_{2}), (x_{1}^{'}, x_{2}^{'})\} \in P_{R(T_{1}) \bigoplus R(T_{2})}$, $\{(x_{1}^{'}, x_{2}^{'}), (y_{1}, y_{2})\} \in (T_{1} \bigoplus T_{2})^{-1}$, $\{(y_{1}, y_{2}), (y_{1}, y_{2})\} \in P_{(N(T_{1}))^{\perp} \bigoplus (N(T_{2}))^{\perp}}.$ Thus, $$\{x, y\} = \{(x_{1}, x_{2}), (y_{1}, y_{2})\} \in (T_{1} \bigoplus T_{2})^{\dagger},$$ which implies $$T_{1}^{\dagger} \bigoplus T_{2}^{\dagger} \subset (T_{1} \bigoplus T_{2})^{\dagger}.$$
To prove the reverse inclusion, let $$\{u, v\} \in (T_{1} \bigoplus T_{2})^{\dagger} = P_{(N(T_{1}))^{\perp} \bigoplus (N(T_{2}))^{\perp}} (T_{1} \bigoplus T_{2})^{-1} P_{R(T_{1}) \bigoplus R(T_{2})},$$ with $u= u^{'} + u^{''},$ where $u^{'} \in R(T_{1}) \bigoplus R(T_{2})$, $u^{''} \in (R(T_{1}) \bigoplus R(T_{2}))^{\perp}.$ Then $$\{u, u^{'}\} \in P_{R(T_{1}) \bigoplus R(T_{2})},~ \{u^{'}, v\} \in (T_{1} \bigoplus T_{2})^{-1},~ \{v, v\} \in P_{(N(T_{1}))^{\perp} \bigoplus (N(T_{2}))^{\perp}}.$$ 
Write $u = (u_{1}, u_{2}) \in H \bigoplus K$, $u^{'} = (u_{1}^{'}, u_{2}^{'})$ with $u_{1} = u_{1}^{'} + u_{1}^{''}$, $u_{2} = u_{2}^{'} + u_{2}^{''},$ where $u_{1}^{'} \in R(T_{1})$, $u_{1}^{''} \in (R(T_{1}))^{\perp}$,  $u_{2}^{'} \in R(T_{2})$, $u_{2}^{''} \in (R(T_{2}))^{\perp}.$\\
Let $v = (v_{1}, {v_{2}}) \in (N(T_{1}))^{\perp} \bigoplus (N(T_{2}))^{\perp}.$\\
This shows that $\{u_{1}^{'}, v_{1}\} \in T_{1}^{-1}$, $\{u_{2}^{'}, v_{2}\} \in T_{2}^{-1}$ which implies $$\{u_{1}, v_{1}\} \in T_{1}^{\dagger}~ \text{and}~ \{u_{2}, v_{2}\} \in T_{2}^{\dagger}.$$ So, $\{u, v\} = \{(u_{1}, u_{2}), (v_{1}, v_{2})\} \in T_{1}^{\dagger} \bigoplus T_{2}^{\dagger},$ and hence, $(T_{1} \bigoplus T_{2})^{\dagger} \subset T_{1}^{\dagger} \bigoplus T_{2}^{\dagger}.$ Therefore, $$(T_{1} \bigoplus T_{2})^{\dagger} = T_{1}^{\dagger} \bigoplus T_{2}^{\dagger}.$$
\end{proof}
\begin{corollary}\label{cor 3.16}
Let $T_{1} \in CR(H_{1}, K_{1})$ and $T_{2} \in CR(H_{2}, K_{2})$ both be closed range linear relations. Then $$\gamma(T_{1} \bigoplus T_{2}) = \min\{\gamma(T_{1}), \gamma(T_{2})\} > 0.$$
\end{corollary}
\begin{proof}
Since, $$(T_{1} \bigoplus T_{2})^{\dagger} = T_{1}^{\dagger} \bigoplus T_{2}^{\dagger},$$ and both $T_{1}^{\dagger}$ and $T_{2}^{\dagger}$ are two bounded operators (by Theorem \ref{thm 3.2}), it follows that $T_{1}^{\dagger} \bigoplus T_{2}^{\dagger}$ is also bounded. Moreover, its norm satisfies $$\|T_{1}^{\dagger} \bigoplus T_{2}^{\dagger}\| = \max\{\|T_{1}^{\dagger}\|, \|T_{2}^{\dagger}\|\}.$$ On the other hand, $R(T_{1} \bigoplus T_{2}) = R(T_{1}) \bigoplus R(T_{2})$ is closed. Hence, by Theorem $\ref{thm 3.4},$ Theorem \ref{thm 3.15}, and Proposition 3.1 \cite{MR3079830}, we obtain 
\begin{align*}
\gamma(T_{1} \bigoplus T_{2}) = \frac{1}{\|T_{1}^{\dagger} \bigoplus T_{2}^{\dagger}\|} & = \frac{1}{\max\{\|T_{1}^{\dagger}\|, \|T_{2}^{\dagger}\|\}}\\ & = \min\{\gamma(T_{1}), \gamma(T_{2})\} > 0.
\end{align*}
\end{proof}
\begin{corollary}\label{cor 3.17}
Let $T_{1} \in CR(H_{1}, K_{1})$ and $T_{2} \in CR(H_{2}, K_{2})$ both be the closed range linear relations. Then $$((T_{1} \bigoplus T_{2})^{\dagger})^{*} = (T_{1}^{\dagger})^{*} \bigoplus (T_{2}^{\dagger})^{*}.$$
\end{corollary}
\begin{proof}
It is easy to show that $$(T_{1} \bigoplus T_{2})^{*} = T_{1}^{*} \bigoplus T_{2}^{*}.$$ By Theorem \ref{thm 3.15}, Theorem \ref{thm 3.1} and Corollary 3.3 (ii)\cite{MR3079830}, we obtain $$((T_{1} \bigoplus T_{2})^{\dagger})^{*} = ((T_{1} \bigoplus T_{2})^{*})^{\dagger} = (T_{1}^{*})^{\dagger} +  (T_{2}^{*})^{\dagger} = (T_{1}^{\dagger})^{*} +  (T_{2}^{\dagger})^{*}.$$
\end{proof}
\begin{theorem}\label{thm 3.18}
Let $T \in CR(H, K)$ be a closed range linear relation. Then
\begin{enumerate}
\item $\vert T^{\dagger} \vert = \vert T^{*}\vert^{\dagger}.$
\item $\vert (T^{\dagger})^{*}\vert = \vert T \vert^{\dagger}.$
\end{enumerate}
\end{theorem}
\begin{proof}
$\it(1)$ Since $T^{\dagger}$ is an operator in domain $K,$ its modulus is well defined and satisfies $$\vert T^{\dagger} \vert = ((T^{\dagger})^{*}T^{\dagger})^{\frac{1}{2}} = ((T^{*})^{\dagger}T^{\dagger})^{\frac{1}{2}} = ((TT^{*})^{\dagger})^{\frac{1}{2}},$$ where the second equality follows from Theorem \ref{thm 3.1} and the last from Theorem \ref{thm 3.3}$\it(6).$ We next establish the general identity $$(S^{\frac{1}{2}})^{\dagger} = (S^{\dagger})^{\frac{1}{2}},$$ where $S$ is a positive self-adjoint linear relation from $K$ into $K$ with closed range. By Proposition 2.4 and Proposition 2.5 \cite{MR3079830}, we have $$N(S^{\frac{1}{2}}) = N(S),~ R(S^{\frac{1}{2}}) = R(S)~ \text{and}~ (N(S))^{\perp} = R(S).$$ Then 
\begin{align*}
((S^{\frac{1}{2}})^{\dagger})^{2} &= (P_{(N(S))^{\perp}} (S^{\frac{1}{2}})^{-1}P_{R(S)})^{2}\\ &= P_{(N(S))^{\perp}}  (S^{\frac{1}{2}})^{-1}  (S^{\frac{1}{2}})^{-1} P_{R(S)}\\ &=  P_{(N(S))^{\perp}} S^{-1}P_{R(S)}\\ &= S^{\dagger}.
\end{align*}
So, $S^{\dagger}$ is non-negative. By Theorem \ref{thm 3.1}, we have $(S^{\dagger})^{*} = (S^{*})^{\dagger} = S^{\dagger}$ is self-adjoint. Theorem 1.5.9 \cite{MR3971207} confirms that $(S^{\dagger})^{\frac{1}{2}} = (S^{\frac{1}{2}})^{\dagger}.$ By Lemma 1.5.8 \cite{MR3971207} and Proposition 2.5 \cite{MR3079830}, we obtain $TT^{*}$ is a non-negative self-adjoint with closed range. Thus, $$\vert T^{*}\vert^{\dagger} = ((TT^{*})^{\frac{1}{2}})^{\dagger} = ((TT^{*})^{\dagger})^{\frac{1}{2}} = \vert T^{\dagger} \vert.$$
$\it(2)$ We consider $T^{*}$ instead of $T$ to get $\vert (T^{\dagger})^{*}\vert = \vert T \vert^{\dagger}.$
\end{proof}
\begin{corollary}\label{cor 3.19}
Let $T_{1} \in CR(H_{1}, K_{1})$ and $T_{2} \in CR(H_{2}, K_{2})$ both be closed range linear relations. Then $$\vert (T_{1} \bigoplus T_{2})^{\dagger}\vert = \vert T_{1}^{\dagger}\vert \bigoplus \vert T_{2}^{\dagger}\vert.$$
\end{corollary}
\begin{proof}
By Theorem \ref{thm 3.15} and Theorem \ref{thm 3.18}, we obtain
\begin{align*}
\vert (T_{1} \bigoplus T_{2})^{\dagger}\vert &= \vert T_{1}^{\dagger} \bigoplus T_{2}^{\dagger}\vert\\
&= ((T_{1}^{\dagger})^{*}T_{1}^{\dagger} \bigoplus (T_{2}^{\dagger})^{*}T_{2}^{\dagger} )^{\frac{1}{2}}\\
&= ((T_{1}T_{1}^{*})^{\dagger} \bigoplus (T_{2}T_{2}^{*})^{\dagger})^{\frac{1}{2}}\\
&= ((T_{1}T_{1}^{*})^{\frac{1}{2}})^{\dagger} \bigoplus ((T_{2}T_{2}^{*})^{\frac{1}{2}})^{\dagger}\\
&= \vert T_{1}^{*}\vert^{\dagger} \bigoplus \vert T_{2}^{*}\vert^{\dagger}\\
&= \vert T_{1}^{\dagger}\vert \bigoplus \vert T_{2}^{\dagger}\vert.
\end{align*}
\end{proof}
\begin{corollary}\label{cor 3.20}
Let $T_{1} \in CR(H_{1}, K_{1})$ and $T_{2} \in CR(H_{2}, K_{2})$ be two closed range linear relations. Then $$\vert T_{1} \bigoplus T_{2} \vert^{\dagger} = \vert T_{1}\vert^{\dagger} \bigoplus \vert T_{2}\vert^{\dagger}.$$
\end{corollary}
\begin{proof}
By Theorem \ref{thm 3.15}, we obtain
\begin{align*}
\vert T_{1} \bigoplus T_{2}\vert^{\dagger} &= ((T_{1}^{*}T_{1} \bigoplus T_{2}^{*}T_{2})^{\frac{1}{2}})^{\dagger} \\
&= ((T_{1}^{*}T_{1})^{\frac{1}{2}})^{\dagger} \bigoplus ((T_{2}^{*}T_{2})^{\frac{1}{2}})^{\dagger}\\
&= \vert T_{1}\vert^{\dagger} \bigoplus \vert T_{2} \vert^{\dagger}.
\end{align*}
\end{proof}
\section{Characterization of the essential spectra of Off-Diagonal Block linear relations in Hilbert Spaces}

In this section, we are concerned with the following essential spectra:
\[
\begin{aligned}
\sigma_{e_1}(T)
&:= \left\{ \lambda \in \mathbb{C} \;:\; T - \lambda I \notin \phi_{+}(H,H) \right\}, \\[6pt]
\sigma_{e_2}(T)
&:= \left\{ \lambda \in \mathbb{C} \;:\; T - \lambda I \notin \phi_{-}(H,H) \right\}, \\[6pt]
\sigma_{e_3}(T)
&:= \left\{ \lambda \in \mathbb{C} \;:\; T - \lambda I \notin \phi_{+}(H, H) \cup \phi_{-}(H, H) \right\}, \\[6pt]
\sigma_{e_4}(T)
&:= \left\{ \lambda \in \mathbb{C} \;:\; T - \lambda I \notin \phi_{+}(H, H) \cap \phi_{-}(H, H) \right\}.
\end{aligned}
\]
\begin{theorem}\label{thm 3.13}
Let $$T = \begin{bmatrix}
0  & A\\
B & 0 \end{bmatrix} = \bigl\{\{(x,y), (y_{a}, x_{b})\}: \{x, x_{b}\} \in B, \{y, y_{a}\} \in A \bigr\}$$ be a linear relation from domain $D(B) \times D(A)$ into $H \times K$, where $A \in CR(K, H)$ and $B \in CR(H, K)$. Assume that $AB$ and $BA$ are both closed linear relations in $H$ and $K,$ respectively. Then $$\rho(T) = \{\lambda \in \mathbb{C}: \lambda^{2} \in \rho(AB) \cap \rho(BA)\},$$ where $\rho(T), \rho(AB), \rho(BA)$ are the resolvent sets of $T, AB$ and $BA$, respectively.
\end{theorem}
\begin{proof}
 We first show that $T$ is closed. Let $\{h, k\} \in \overline{T}.$ Then there exist a sequence $$\{(x_{n}, y_{n}), (y_{an}, x_{bn})\} \to \{h, k\}~ \text{as}~ n \to \infty,$$ where $\{x_{n}, x_{bn}\} \in B$ and $\{y_{n}, y_{an}\} \in A$, for all $n \in \mathbb{N}.$ So, $\{x_{n}\}$, $\{y_{n}\}$, $\{y_{an}\}$ and $\{x_{bn}\}$ all are Cauchy sequences, there exist $x, y, y_{a}, x_{b}$ such that $$\{x_{n}\} \to x,~ \{y_{n}\} \to y,~ \{y_{an}\} \to y_{a}~ \text{and}~ \{x_{bn}\} \to x_{b}~ \text{as}~ n \to \infty.$$ Since $A$ and $B$ both are closed, then $h = (x, y)$ and $k = (y_{a}, x_{b})$, where $\{x, x_{b}\} \in B$, $\{y, y_{a}\} \in A.$ It confirms that $\{h, k\} \in T.$ Thus, $T$ is closed.
 \medskip
\noindent
 Now, we prove that $\rho(T) \subset \{\lambda \in \mathbb{C}: \lambda^{2} \in \rho(AB) \cap \rho(BA)\}.$\\
$Case(1):$ Let $\lambda \in \rho(T) \setminus \{0\}.$ Then $(T- \lambda)^{-1}$ is bounded operator in domain $H \times K.$ It is easy to show that $$(T - \lambda) = \begin{bmatrix} 
-\lambda I_{H} & A\\
B  & -\lambda I_{K}\end{bmatrix},$$ where $I_{H}$ and $I_{K}$ are two identity operators in $H$ and $K$, respectively.\\
Let $(g_{1}, g_{2}) \in H \times K.$ Then there exists $(f_{1}, f_{2})$ in $D(T)$ such that $$\{(f_{1}, f_{2}), (g_{1}, g_{2})\} \in (T - \lambda),$$ where $g_{1} = -\lambda f_{1} + f_{2a}$, $g_{2} = -\lambda f_{2} + f_{1b}$ with $\{f_{1}, f_{1b}\} \in B$, $\{f_{2}, f_{2a}\} \in A.$\\
Thus, $f_{2a} = g_{1} + \lambda f_{1}$, $f_{1b} = g_{2} + \lambda f_{2}.$ Taking $g_{1} = 0$, then for $(0, g_{2})$, we obtain $$\{(f_{1}^{'}, f_{2}^{'}), (0, g_{2})\} \in (T - \lambda)$$ such that $$f_{2a}^{'} = \lambda f_{1}^{'},~~ f_{1b}^{'} = g_{2} + \lambda f_{2}^{'},$$ where $\{f_{1}^{'}, f_{1b}^{'}\} \in B$ and $\{f_{2}^{'}, f_{2a}^{'}\} \in A.$\\
Now, we claim that $\{f_{2}^{'}, g_{2}\} \in \frac{1}{\lambda}(BA - \lambda^{2}).$ This is true that $$\{\lambda f_{1}^{'}, \lambda f_{1b}^{'}\} = \{f_{2a}^{'}, \lambda g_{2} + \lambda^{2} f_{2}^{'}\} \in B$$ which implies $$\{f_{2}^{'}, \lambda g_{2} + \lambda^{2}f_{2}^{'}\} \in BA.$$ Thus, $\{f_{2}^{'}, g_{2}\} \in \frac{1}{\lambda}(BA - \lambda^{2}).$ So, $D((BA - \lambda^{2})^{-1}) = K.$\\
Again, let $\{w, 0\} \in \frac{1}{\lambda}(BA - \lambda^{2}).$ Then $\{w, \lambda^{2}w\} \in BA$. There exists $s \in R(A)$ such that $$\{w, s\} \in A~ \text{and}~ \{s, \lambda^{2}w\} \in B.$$ Hence, $\{(\frac{1}{\lambda} s, w), (0, 0)\} \in T- \lambda.$ $(T - \lambda)$ is one-one because $(T- \lambda)^{-1}$ is an operator. It proves that $s = w= 0,$ which implies $(BA - \lambda^{2})^{-1}$ is an operator.\\ Since $BA$ is closed, then $(BA - \lambda^{2})^{-1}$ is also closed in domain $K$. By the closed graph theorem (Theorem III.4.2 \cite{MR1631548}), we have $\lambda^{2} \in \rho(BA).$\\ Similarly, it is evident that $\lambda^{2} \in \rho(AB)$. Hence, $$\lambda^{2} \in \rho(AB) \cap \rho(BA).$$
$Case(2):$ Let $0 \in \rho(T).$ Then $T^{-1}$ is a bounded operator in domain $H \times K.$ For an element $(s_{1}, s_{2}) \in H \times K$, we get an unique $(r_{1}, r_{2}) \in H \times K$ such that $$(s_{1}, s_{2}) = (r_{2a}, r_{1b}), ~\text{where}~ \{r_{2}, r_{2a}\} \in A, \{r_{1}, r_{1b}\} \in B.$$
So, $D(A^{-1}) = H$ and $D(B^{-1}) = K$ with $A^{-1}$ and $B^{-1}$ is closed because the closeness of $A$ and $B$. Moreover, $A^{-1}$ and $B^{-1}$ both are operators because of the operator $$T^{-1} = \begin{bmatrix}
0 & B^{-1} \\
A^{-1} & 0\end{bmatrix}.$$ Thus, $A^{-1}$ and $B^{-1}$ are two bounded operators with domains $H$ and $K$, respectively. Again, $(AB)^{-1} = B^{-1}A^{-1}$ and $(BA)^{-1} = A^{-1}B^{-1}$ are two  bounded operators in domain $H$ to $H$ and $K$ to $K$ respectively. It shows that $0 \in \rho(AB) \cap \rho(BA).$ Hence,
\begin{align}
\lambda \in \rho(T) \implies \lambda^{2} \in \rho(AB) \cap \rho(BA).
\end{align}
Conversely, Suppose that $\mu^{2} \in \rho(AB) \cap \rho(BA).$ Then $$(AB - \mu^{2})^{-1},~ \text{and}~ (BA- \mu^{2})^{-1}$$ both are bounded operators in the domain $H$ into $H$ and the domain $K$ into $K$, respectively.\\
Let $\{x, y\} \in T^{2}.$ Then there exists $z \in R(T)$ such that $$\{x, z\} \in T,~ \text{and} \{z, y\} \in T.$$ Taking $\{x, z\} = \{(u, v), (v_{a}, u_{b})\}$ and $\{z, y\} = \{(v_{a}, u_{b}), (u_{ba}, v_{ab})\}$, with $$\{u, u_{b}\} \in B,~ \{v, v_{a}\} \in A,~  \{v_{a}, v_{ab}\} \in B,~ \{u_{b}, u_{ba}\} \in A.$$
Thus, $\{x, y\} = \{(u, v), (u_{ba}, v_{ab})\} \in \begin{bmatrix}
AB & 0\\
0 & BA \end{bmatrix}.$\\
Again, consider $\{p, q\} \in \begin{bmatrix}
AB & 0\\
0 & BA \end{bmatrix}.$ Write $\{p, q\} = \{(w, s), (w_{ab}, s_{ba})\}$ with $\{w, w_{ab}\} \in AB$ and $\{s, s_{ba}\} \in BA,$ then there exist $w_{b}$ and $s_{a}$ such that $$\{w, w_{b}\} \in B,~ \{w_{b}, w_{ab}\} \in A,~ \{s, s_{a}\} \in A,~ \text{and}~ \{s_{a}, s_{ba}\} \in B.$$ 
This implies that $\{(w, s), (s_{a}, w_{b})\} \in T$ and $\{(s_{a}, w_{b}), (w_{ab}, s_{ba})\} \in T$. So, $\{p, q\} \in T^{2}.$ Since $AB$ and $BA$ both are closed, then $$T^{2} = \begin{bmatrix}
AB & 0\\
0 & BA \end{bmatrix}$$ is closed. Furthermore, $T^{2} - \mu^{2} = \begin{bmatrix}
AB- \mu^{2} & 0\\
0 & BA - \mu^{2} \end{bmatrix}$ implies $\mu^{2} \in \rho(T^{2}).$\\ Now, we will show that $T^{2} - \mu^{2} = (T + \mu) (T - \mu) = (T - \mu )(T+ \mu).$\\
Let $\{u^{'}, v^{'}\} \in (T + \mu) (T - \mu).$ Then there exists $s^{'} \in R(T - \mu)$ such that $$\{u^{'}, s^{'}\} \in T- \mu,~ \{s^{'}, v^{'}\} \in T + \mu.$$
The condition $\{u^{'}, s^{'}\} \in T- \mu$ confirms that $$\{ u^{'},  s^{'} + \mu u^{'}\} \in T~ \text{and}~ \{\mu u^{'}, \mu s^{'} + \mu^{2}u^{'}\} \in T,$$ and the condition $\{s^{'}, v^{'}\} \in T + \mu$ shows that $$\{s^{'}, v^{'} - \mu s^{'}\} \in T.$$
So, $\{\mu u^{'} + s^{'}, v^{'} + \mu^{2}u^{'}\} \in T$ which implies $\{u^{'}, v^{'} + \mu^{2}u^{'}\} \in T^{2}.$ It proves that $\{u^{'}, v^{'}\} \in T^{2} - \mu^{2}.$ Hence, $$(T + \mu) (T - \mu) \subset T^{2} - \mu^{2}.$$
Let $\{u^{''}, v^{''}\} \in T^{2} - \mu^{2}.$ Then $\{u^{''}, v^{''} + \mu^{2}u^{''}\} \in T^{2}.$ There exists $u^{'''} \in R(T)$ such that $$\{u^{''}, u^{'''}\} \in T~ \text{and}~ \{u^{'''}, v^{''} + \mu^{2}u^{''}\} \in T.$$ Furthermore, $\{u^{''}, u^{'''} - \mu u^{''}\} \in T- \mu$ and $\{u^{'''} - \mu u^{''}, v^{''} + \mu^{2}u^{''}- \mu u^{'''}\} \in T$. So, $\{u^{'''} - \mu u^{''}, v^{''}\} \in T + \mu.$ It confirms that $\{u^{''}, v^{''}\} \in (T + \mu) (T - \mu).$ We are ready to say that $$T^{2} - \mu^{2} = (T + \mu) (T - \mu).$$ 
Similarly, it can be proven that $T^{2} - \mu^{2} = (T - \mu) (T + \mu).$  By the property, $$\mu^{2} \in \rho(T^{2})~ \text{and}~  T^{2} - \mu^{2} = (T + \mu) (T - \mu) = (T - \mu) (T + \mu),$$ we have $N(T - \mu) = \{0\}$, and  $R(T - \mu) = H \times K.$ Since $T$ is closed, then  $(T - \mu)^{-1}$ is a closed. Thus, $(T-\mu)^{-1}$ is a bounded operator in domain $H \times K,$ and hence $\mu \in \rho(T).$ Therefore, $\rho(T) = \{\lambda \in \mathbb{C}: \lambda^{2} \in \rho(AB) \cap \rho(BA)\}.$
\end{proof}
\begin{corollary}\label{cor 3.14}
Let $T \in CR(H, K) \cap BR(H, K)$ be a closed range linear relation. Define $$\mathcal{A} = \begin{bmatrix} 
0 & T^{\dagger}\\
T & 0 \end{bmatrix}.$$ Then the resolvent set of $A$ is given by $$\rho(\mathcal{A}) = \{\lambda \in \mathbb{C}: \lambda^{2} \in \rho(TT^{\dagger})  \cap \rho(T^{\dagger}T)\}.$$
\end{corollary}
\begin{proof}
By Theorem \ref{thm 3.6}, we have $TT^{\dagger}$ and $T^{\dagger}T$ both are closed. Using Theorem \ref{thm 3.13}, we get $\rho(\mathcal{A}) = \{\lambda \in \mathbb{C}: \lambda^{2} \in \rho(TT^{\dagger})  \cap \rho(T^{\dagger}T)\}.$
\end{proof}
By Proposition 1.3.9(ii) \cite{MR3971207}, it follows that $(TS)^{*} = S^{*}T^{*},$ when $S \in LR(H, K)$ and $T$ is a bounded operator from $K$ into $H.$ The next lemma shows that the identity $(TS)^{*} = S^{*}T^{*}$ is also valid when $S \in BCR(H, H)$ and $T \in CR(H, H)$ with $0 \in \rho(S).$
\begin{lemma}\label{lemma 4.3}
Let $T \in CR(H, H)$ and $S \in BCR(H, H).$ Assume $0 \in \rho(S).$ Then 
$$(TS)^{*} = S^{*}T^{*}.$$
\end{lemma}
\begin{proof}
It follows from Proposition 1.3.9(ii) \cite{MR3971207} that $S^{*}T^{*} \subset (TS)^{*}.$ To prove the reverse inclusion, let $\{x, y\} \in T.$ Then $x \in H = D(S^{-1}).$ There exists $z \in H$ such that $\{z, x\} \in S$ which implies $\{z, y\} \in TS.$ So, $\{x, y\} \in TSS^{-1}.$ This shows that $$T \subset TSS^{-1}.$$
Taking adjoint both sides, we obtain
$$(S^{-1})^{*}(TS)^{*} \subset (TSS^{-1})^{*} \subset T^{*}.$$
It confirms that 
$$(S^{-1}S)^{*}(TS)^{*} = S^{*}(S^{-1})^{*} (TS)^{*} \subset S^{*}T^{*}.$$
Using Theorem \ref{thm 3.6}, we get $S^{-1}S$ is a bounded operator in domain $H,$ hence\\ $S^{-1}S = I_{H}.$ Thus, $(S^{-1}S)^{*} = I_{H}.$ Therefore, $$(TS)^{*} \subset S^{*}T^{*}.$$
\end{proof}
\begin{theorem}\label{thm 4.4}
Let $T \in CR(H, H)$ and $U$ be an unitary operator in $H.$ Then the following statements hold:
\begin{enumerate}
\item $T \in \phi_{+}(H, H)$ ~\text{if and only if}~ $U^{*}TU \in \phi_{+}(H,H).$
\item $T \in \phi_{-}(H, H)$ ~\text{if and only if}~ $U^{*}TU \in \phi_{-}(H,H).$
\item $T \in \phi_{+}(H, H) \cup \phi_{-}(H, H)$ ~\text{if and only if}~ $U^{*}TU \in \phi_{+}(H,H) \cup \phi_{-}(H, H).$
\item $T \in \phi_{+}(H, H) \cap \phi_{-}(H, H)$ ~\text{if and only if}~ $U^{*}TU \in \phi_{+}(H,H) \cap \phi_{-}(H, H).$
\end{enumerate}
\end{theorem}
\begin{proof}
$\it(1)$ Suppose that $ T \in \phi_{+}(H, H).$ Then $R(T)$ is closed and $\dim(N(T)) < \infty.$ It is evident that $N(U^{*}TU) = N(TU).$\\
From Lemma 3.3 (iii)\cite{MR00027}, it follows that $T^{*} \in \phi_{-}(H, H).$ Again, $U^{*} \in \phi_{-}(H, H),$ since  $(R(U^{*}))^{\perp} =  \{0\}.$ By Lemma \ref{lemma 4.3}, we obtain $(TU)^{*} = U^{*}T^{*},$ hence $U^{*}T^{*}$ is closed. Using Proposition 6.1(ii) \cite{MR00027}, we get $(TU)^{*} \in \phi_{-}(H, H)$ which implies that $TU \in \phi_{+}(H, H).$ Then $\dim(N(U^{*}TU)) = \dim(N(TU)) < \infty$ and $R(TU)$ is closed. Now, we claim that $R(U^{*}TU)$ is closed.\\
Let $y \in \overline{R(U^{*}TU)}.$ Then there exists $\{x_{n}, y_{n}\} \in U^{*}TU$ such that $y_{n} \to y,$ as $n \to \infty.$ 
Again we get $z_{n} \in R(TU)$ such that $\{x_{n}, z_{n}\} \in TU,$ and $\{z_{n}, y_{n}\} \in U^{*}$ for all $n \in \mathbb{N}.$\\
From the identity $\|z_{n}- z_{m}\| = \|y_{n}-y_{m}\|$ for all $n,m \in \mathbb{N},$ we have that $\{z_{n}\}$ is Cauchy sequence. There exists $z \in H$ such that $z_{n} \to z$ as $n \to \infty.$ Since, $R(TU)$ is closed, hence $z \in R(TU).$ It shows that $\{x, z\} \in TU$ and $\{z, y\} \in U^{*},$ for some $x \in H.$\\
Thus $y \in R(U^{*}TU)$ which implies $R(U^{*}TU)$ is closed. 
Therefore, $U^{*}TU \in \phi_{+}(H, H).$\\
The converse is obvious because $T = UU^{*}TUU^{*}.$\\

$\it(2)$ Using Lemma \ref{lemma 4.3} and Lemma 3.3 \cite{MR00027}, we obtain
$T \in \phi_{-}(H,H)$ if and only if $T^{*} \in \phi_{+}(H,H)$ is and only if $U^{*}T^{*}U \in \phi_{+}(H,H)$ if and only if $U^{*}TU \in \phi_{-}(H,H).$\\

$\it(3)$ and $\it(4)$ both are true from the above two statements.
\end{proof}
\begin{theorem}\label{thm 4.5}
Let $$T = \begin{bmatrix}
0  & A\\
B & 0 \end{bmatrix} = \bigl\{\{(x,y), (y_{a}, x_{b})\}: \{x, x_{b}\} \in B, \{y, y_{a}\} \in A \bigr\}$$ be a linear relation from domain $H \times K,$ where $A \in CR(K, H) \cap BCR(K, H) $ and $B \in CR(H, K) \cap BCR(H, K).$ Assume that $AB$ and $BA$ are both closed linear relations in $H$ and $K$ respectively. Then
$$\sigma_{e_{1}}(T) = \{\lambda \in \mathbb{C}: \lambda^{2} \in \sigma_{e_{1}}(AB) \cup \sigma_{e_{1}}(BA) \}.$$
\end{theorem}
\begin{proof}
It is evident to show that 
$$-T - \gamma =  \begin{bmatrix}
-\gamma  & -A\\
-B & -\gamma \end{bmatrix} =  \begin{bmatrix}
I & 0\\
0 & -I \end{bmatrix}  \begin{bmatrix}
-\gamma  & A\\
B & -\gamma \end{bmatrix} \begin{bmatrix}
I & 0\\
0 & -I \end{bmatrix}.$$
Here, $\begin{bmatrix}
I & 0\\
0 & -I \end{bmatrix}$ is an unitary operator. By Theorem \ref{thm 4.4} we obtain $ \gamma \notin \sigma_{e_{1}}(T)$ if and only if $ \gamma \notin \sigma_{e_{1}}(-T).$\\
Let $ \lambda \notin \sigma_{e_{1}}(T).$ Then $(T - \lambda)$ and $(T + \lambda)$ both are in $\phi_{+}(H \times K, H \times K).$ Using Proposition 6.1(i) \cite{MR00027}, we get $(T^{2} - \lambda^{2}) \in \phi_{+}(H \times K, H \times K)$ which implies $(AB - \lambda^{2})$  and $(BA - \lambda^{2})$ both are in $\phi_{+}(H, H)$ and $\phi_{+}(K, K)$ respectively. Thus, $$\lambda^{2} \notin \sigma_{e_{1}}(AB) \cup \sigma_{e_{1}}(BA).$$ It shows that $\{\lambda \in \mathbb{C}: \lambda^{2} \in \sigma_{e_{1}}(AB) \cup \sigma_{e_{1}}(BA) \} \subset \sigma_{e_{1}}(T).$\\
To prove the converse inclusion, let $\mu^{2} \notin \sigma_{e_{1}}(AB) \cup \sigma_{e_{1}}(BA).$ Then $$(T^{2} - \mu^{2}) \in \phi_{+}(H \times K, H \times K).$$ Using Proposition 6.3 \cite{MR00027}, we get $(T - \mu) \in \phi_{+}(H \times K, H \times K)$ which implies that $\mu \notin \sigma_{e_{1}} (T).$ Therefore, $$\sigma_{e_{1}}(T) = \{\lambda \in \mathbb{C}: \lambda^{2} \in \sigma_{e_{1}}(AB) \cup \sigma_{e_{1}}(BA) \}.$$
\end{proof}
\begin{theorem}\label{thm 4.6}
Let $$T = \begin{bmatrix}
0  & A\\
B & 0 \end{bmatrix} = \bigl\{\{(x,y), (y_{a}, x_{b})\}: \{x, x_{b}\} \in B, \{y, y_{a}\} \in A \bigr\}$$ be a linear relation from domain $H \times K,$ where $A \in CR(K, H) \cap BCR(K, H) $ and $B \in CR(H, K) \cap BCR(H, K).$ Assume that $AB$ and $BA$ are both closed linear relations in $H$ and $K$ respectively. Then
$$\sigma_{e_{2}}(T) = \{\lambda \in \mathbb{C}: \lambda^{2} \in \sigma_{e_{2}}(AB) \cup \sigma_{e_{2}}(BA) \}.$$
\end{theorem}
\begin{proof}
Let $\lambda \notin \sigma_{e_{2}}(T).$ Using Theorem \ref{thm 4.4}(ii), we obtain $\lambda \notin \sigma_{e_{2}}(-T)$ which implies that $(T^{2} - \lambda^{2}) \in \phi_{-} (H \times K, H \times K).$ This shows that
$$\lambda^{2} \notin \sigma_{e_{2}}(AB) \cup \sigma_{e_{2}}(BA).$$ Hence,
 $\{\lambda \in \mathbb{C}: \lambda^{2} \in \sigma_{e_{2}}(AB) \cup \sigma_{e_{2}}(BA) \} \subset \sigma_{e_{2}}(T).$\\
 To prove the reverse inclusion, consider $\mu^{2} \notin \sigma_{e_{2}}(AB) \cup \sigma_{e_{2}}(BA).$ Then $(AB - \mu^{2})$ and $(BA - \mu^{2})$ both are closed range linear relations with finite codimensions.\\ It shows that $(T- \mu) (T + \mu) = (T^{2} - \mu^{2}) \in \phi_{-}(H \times K, H \times K).$ By Proposition 6.4, we get $(T - \mu) \in \phi_{-}(H \times K, H \times K).$ It confirms that $\mu \notin \sigma_{e_{2}}(H \times K).$ Therefore,
 $$\sigma_{e_{2}}(T) = \{\lambda \in \mathbb{C}: \lambda^{2} \in \sigma_{e_{2}}(AB) \cup \sigma_{e_{2}}(BA) \}.$$
 \end{proof}
 \begin{corollary}\label{cor 4.7}
 Let $$T = \begin{bmatrix}
0  & A\\
B & 0 \end{bmatrix} = \bigl\{\{(x,y), (y_{a}, x_{b})\}: \{x, x_{b}\} \in B, \{y, y_{a}\} \in A \bigr\}$$ be a linear relation from domain $H \times K,$ where $A \in CR(K, H) \cap BCR(K, H) $ and $B \in CR(H, K) \cap BCR(H, K).$ Assume that $AB$ and $BA$ are both closed linear relations in $H$ and $K$ respectively. Then the following statements hold:
\begin{enumerate}
\item $\sigma_{e_{4}}(T) = \Bigl \{\lambda \in \mathbb{C}: \lambda^{2} \in \sigma_{e_{4}}(AB) \cup \sigma_{e_{4}}(BA) \Bigr\}.$
\item $\sigma_{e_3}(T)
=
 \Bigl \{
\lambda \in \mathbb{C} :
\lambda^{2} \in
\sigma_{e_3}(AB) \cup \sigma_{e_3}(BA) \cup \bigl[\sigma_{e_1}(AB) \cap \sigma_{e_2}(BA)\bigr]\\
\cup \bigl[\sigma_{e_2}(AB) \cap \sigma_{e_1}(BA)\bigr]
\Bigr\}.
$
\end{enumerate}
\end{corollary}
\begin{proof}
    $\it(1)$ The proof directly follows from Theorem \ref{thm 4.5} and Theorem \ref{thm 4.6}.\\
   
     $\it(2)$ \begin{align*}
    &\sigma_{e_{3}}(T)\\ =& \sigma_{e_{1}}(T) \cap \sigma_{e_{2}}(T)\\ = & \{\lambda \in \mathbb{C}: \lambda^{2} \in \sigma_{e_{1}}(AB) \cup \sigma_{e_{1}}(BA) \}\cap \{\lambda \in \mathbb{C}: \lambda^{2} \in \sigma_{e_{2}}(AB) \cup \sigma_{e_{2}}(BA) \}\\ =& \left\{
\lambda \in \mathbb{C} :
\lambda^{2} \in
\sigma_{e_3}(AB)
\cup \sigma_{e_3}(BA)
\cup \bigl[\sigma_{e_1}(AB) \cap \sigma_{e_2}(BA)\bigr]
\cup \bigl[\sigma_{e_2}(AB) \cap \sigma_{e_1}(BA)\bigr]
\right\}.
    \end{align*}
\end{proof}
\begin{corollary}\label{cor 4.8}
Let $T \in CR(H, K)$ be a closed range linear relation with $D(T) = H$. Define $$\mathcal{A} = \begin{bmatrix} 
0 & T^{\dagger}\\
T & 0 \end{bmatrix}.$$ Then the following statements hold:
\begin{enumerate}
\item $\sigma_{e_{1}}(\mathcal{A}) = \Bigl\{\lambda \in \mathbb{C}: \lambda^{2} \in \sigma_{e_{1}}(TT^{\dagger}) \cup \sigma_{e_{1}}(T^{\dagger}T) \Bigr \}.$
\item $\sigma_{e_{2}}(\mathcal{A}) = \Bigl\{\lambda \in \mathbb{C}: \lambda^{2} \in \sigma_{e_{2}}(TT^{\dagger}) \cup \sigma_{e_{2}}(T^{\dagger}T) \Bigr\}.$
\item $\sigma_{e_3}(\mathcal{A})=
 \Bigl\{ \lambda \in \mathbb{C} : \lambda^{2} \in \sigma_{e_3}(TT^{\dagger}) \cup \sigma_{e_3}(T^{\dagger}T)  \cup \bigl[\sigma_{e_1}(T^{\dagger}T) \cap \sigma_{e_2}(TT^{\dagger})\bigr] \\
\cup \bigl[\sigma_{e_2}(T^{\dagger}T) \cap \sigma_{e_1}(TT^{\dagger})\bigr]
\Bigr\}.$
\item $\sigma_{e_{4}}(\mathcal{A}) = \Bigl \{\lambda \in \mathbb{C}: \lambda^{2} \in \sigma_{e_{4}}(TT^{\dagger}) \cup \sigma_{e_{4}}(T^{\dagger}T) \Bigr\}.$
\end{enumerate}
\end{corollary}
\begin{proof}
By Theorem \ref{thm 3.6}, we have $TT^{\dagger}$ and $T^{\dagger}$ both are closed. Moreover, $T \in CR(H,K)$ and $T^{\dagger}$ is a bounded operator from $K$ to $H.$\\ Using Theorem \ref{thm 4.5}, Theorem \ref{thm 4.6}, and Corollary \ref{cor 4.7}, we obtain the validity of statements $\it(1),$ $\it(2),$ $\it(3),$ and $\it(4).$
\end{proof}

\section{Conclusions}
 In this paper, we have characterized the essential spectra and the resolvent set of the off-diagonal block linear relation $\begin{bmatrix}
0 & A \\
B & 0
\end{bmatrix}$
in terms of the essential spectra and the resolvent sets of the linear relations $AB$ and $BA$. The obtained results establish a clear and precise relationship between the spectral properties of the block linear relation and those of the associated compositions.

We have further investigated the Moore--Penrose inverses of closed linear relations in Hilbert spaces and applied these results to the spectral analysis of the off-diagonal block linear relation $\mathcal{A} =
\begin{bmatrix}
0 & T^{\dagger} \\
T & 0
\end{bmatrix},$
where $T$ is a closed, continuous linear relation with closed range between Hilbert spaces and $T^{\dagger}$ denotes its Moore--Penrose inverse.\\ Furthermore, by proving the result $(T_{1} \bigoplus T_{2})^{\dagger} = T_{1}^{\dagger} \bigoplus T_{2}^{\dagger}$, the Moore-Penrose inverses of the direct sum of closed linear relations $T_{1} \text{ and } T_{2}$ with closed ranges in Hilbert spaces, we obtain $\gamma(T_{1} \bigoplus T_{2}) = \min \{\gamma(T_{1}),\gamma(T_{2})\} > 0$ and a few more characterizations.
These findings contribute to a deeper understanding of the interplay between block linear relation structures, the composition of two linear relations, and generalized inverses, and they provide a foundation for further investigations of spectral properties of structured linear relations in Hilbert spaces.

\section*{Author Contributions} 
This is a single-authored publication which was solely prepared by the author.
\section*{Funding}
No funding was received to assist with the preparation of this manuscript.
\section*{Data Availability}
Not applicable

\section*{Declarations}
  The author declares that there are no conflicts of interest.

\end{document}